# Well-posed variational formulations of Friedrichs-type systems

Martin Berggren, Linus Hägg

*Department of Computing Science, Umeå University, Sweden*

May 3, 2021

**Abstract**

All finite element methods, as well as much of the Hilbert-space theory for partial differential equations, rely on variational formulations, that is, problems of the type: find $u \in V$ such that $a(v, u) = l(v)$ for each $v \in L$, where $V$, $L$ are Sobolev spaces. However, for systems of Friedrichs type, there is a sharp disparity between established well-posedness theories, which are not variational, and the very successful discontinuous Galerkin methods that have been developed for such systems, which are variational. In an attempt to override this dichotomy, we present, through three specific examples of increasing complexity, well-posed variational formulations of boundary and initial–boundary-value problems of Friedrichs type. The variational forms we introduce are generalizations of those used for discontinuous Galerkin methods, in the sense that inhomogeneous boundary and initial conditions are enforced weakly through integrals in the variational forms. In the variational forms we introduce, the solution space is defined as a subspace $V$ of the graph space associated with the differential operator in question, whereas the test function space $L$ is a tuple of $L^2$ spaces that separately enforce the equation, boundary conditions of characteristic type, and initial conditions.

## 1 Introduction

Many mathematical models in applications are most naturally derived and formulated as systems of first-order partial differential equations, for instance the Maxwell equations and the linearized Euler equations of gas dynamics. To analyze broad classes of such systems, Friedrichs [11] introduced the concept of symmetric, positive systems. One attractive aspect of this concept is that it defies the type classification of equations as elliptic, parabolic, or hyperbolic. Indeed, the motivation for Friedrichs to develop his approach was to find a unified framework that encompasses equations that change type, such as the equations of transonic flow. Friedrichs' approach has been developed and extended during the years, for instance by Lax & Phillips [15] and Rauch [21].

More recently, as demonstrated for instance by Houston et al. [12] and by the extensive investigations in Jensen's Ph. D. thesis [13], there have been a renewed interest in the theory of Friedrichs systems, due to the development of discontinuous Galerkin methods, which have emerged as particularly suitable numerical methods for systems written in first-order form.

Of particular relevance for our contribution is the reformulation and abstraction of Friedrichs systems due to Ern, Guermond, and Caplain (EGC) [8], which we now briefly review. Assume we want to solve a system of differential equations

$$Tu = f, \qquad (1.1)$$

supplemented with suitable side conditions. EGC consider a real Hilbert space $L$, equipped with inner product $(\cdot, \cdot)_L$ and norm $\|\cdot\|_L$, and a dense subspace $\mathscr{D}$ of $L$. Typically, $L$ will be an $L^2$ space and $\mathscr{D}$ smooth functions with compact support. In this theory, operator $T$ is assumed to satisfy the bound

$$\|(T + \widetilde{T})\phi\|_L \leq C \|\phi\|_L \qquad \forall \phi \in \mathscr{D}, \qquad (1.2)$$

---

where $\widetilde{T}$ is the formal adjoint of $T$ with respect to $L$, that is, the operator satisfying

$$(T\phi, \psi)_L = (\phi, \widetilde{T}\psi)_L \qquad \forall \phi, \psi \in \mathscr{D}. \tag{1.3}$$

We may always write $T = (T - \widetilde{T})/2 + (T + \widetilde{T})/2$; that is, the operator can always be viewed as a sum of a formally skew symmetric and a formally symmetric operator. Thus, condition (1.2) means that the operators that are considered within this theory are those whose symmetric part is bounded in $L$. No such requirement is assumed on the other part; in a sense, the operator is required to be "essentially" formally skew-symmetric. An operator that satisfies condition (1.2) for $L = L^2(0, 1)$ is $1 + \partial_x$, and an operator that does not is $1 - \partial_{xx}$. Thus, condition (1.2) is tailored for equations in first-order form.

Solutions to equation (1.1) are considered in a subspace $V$ of the graph space

$$W = \{u \in L \mid Tu \in L\}, \tag{1.4}$$

chosen so that $T : V \to L$ is an isomorphism. For boundary-value problems, the space $V$ is directly tied to the choice of boundary conditions. EGC introduce an abstract characterization of these, which is then verified on a case by case basis as the theory is applied to various boundary-value problems. When this framework, as well as the classical Friedrichs theory [21], is applied to specific boundary- or initial–boundary-value problem, it leads to formulations involving what in the finite-element community are known as *essential*, or *strongly enforced* boundary and initial conditions. That is, the boundary conditions are build into the definition of the solution space $V$. Inhomogeneous conditions then need to be treated by a lifting procedure, so that they will be incorporated in the right-hand side $f$.

A related but more comprehensive framework, in which skewness also is central, is the extensive Hilbert space solution theory [19, Chapter 6], originating from the work of Picard [18], which considers so-called evo-systems of the form

$$\partial_t v + Au = f, \tag{1.5a}$$
$$v = Mu, \tag{1.5b}$$

where $\partial_t$ denotes the time derivative, $A$ typically is a linear skew-selfadjoint spatial differential operator, and $M$ is a bounded linear operator. One of the merits of the theory for evo-systems is that it incorporates broad classes of material laws (1.5b), and it allows generalizations to quite complicated initial–boundary-value problems [20]. Moreover, similarly as the original theory for Friedrichs systems, this theory can encompass equations of changing type, a property used by Franz et al. [10] to propose a numerical method for an unsteady equation of changing type.

The theories outlined above are not *variational* in the sense of the standard Hilbert-space theory for partial differential equations. In the variational framework, extensively used, for instance, in the book by Wloka [26], a linear boundary-value problem is reformulated into a problem of the following kind:

$$\begin{aligned} &\text{Find } u \in V \text{ such that} \\ &a(v, u) = l(v) \qquad \forall v \in L, \end{aligned} \tag{1.6}$$

where $a$ is a continuous bilinear form defined on a pair of suitable Sobolev spaces $V$ and $L$, and $l$ is a bounded linear functional on $L$. The well-posedness of problem (1.6) is characterized by the conditions of theorem 2.1 below. In the common case of $L = V$, sufficient conditions are given by the celebrated Lax–Milgram lemma.

*Remark* 1.1. The historical use of the term "variational" comes from the calculus of variations, in which the variational problem constitutes the Euler–Lagrange equations of stationarity of a functional. Here we use the term in a wider sense, for all equations of the type (1.6), whether they are Euler–Lagrange equations or not.

*Remark* 1.2. Note that we here position the test function as the *first* argument of the bilinear form, a convention not shared by all authors.

The variational approach can also be generalized to time dependent problems. J. L. Lions generalized the Lax–Milgram lemma to make it applicable to parabolic initial–boundary-value problems [25,

Lemma 41.2]. Otherwise, in combination with Galerkin approximations or by using the semi-group approach, variational methods can also be used to provide well-posedness results for parabolic as well as for second-order hyperbolic problems [5, 26].

However, variational forms play a much more limited role in the classical theory for first-order partial differential equations. A variational form involving the formal adjoint operator and a so-called semi-admissible boundary operator is used by Jensen [13, Thm. 28] to establish existence of weak solutions to Friedrichs systems. A similar procedure, but with strong enforcement of the semi-admissible boundary conditions, was used already in the original contribution by Friedrichs [11]. Evans [9, § 7.3] also employs a variational form in order to prove existence of a weak solution to the pure initial-value problem for a symmetric hyperbolic system. However, in contrast to the variational theory for second-order problems, uniqueness and continuous dependency of data do not follow from these variational formulations. These remaining aspects of well-posedness of the problem need to be demonstrated separately.

The finite-element method always relies on variational techniques, which is why it was initially developed for partial differential equations that naturally are analyzed as variational problems, such as the equations of linear elasticity. The discontinuous Galerkin method for Friedrichs-type systems is therefore unusual, as it introduces a variational problem in the discrete setting, whose solution converges, when the discretization is refined, to an exact solution that is not constructed using strictly variational means. The main idea behind discontinuous Galerkin methods is to impose boundary conditions, homogeneous or inhomogeneous, as well as interelement continuity weakly through integrals added to the variational form. In contrast, in Friedrichs-type solution theory, boundary conditions, by necessity homogeneous, are typically imposed in the definition of the solution space. In their article on discontinuous Galerkin methods for Friedrichs systems, Ern & Guermond [7, § 2.3] show that the solution to the original problem, *constructed by their non-variational approach* nevertheless uniquely solves a variational problem of the type (1.6) with $L = V = W$ and for $l$ given by $l(v) = (f, v)_L$. However, as demonstrated below in § 3.2, the variational problem, as given by the authors, will be ill posed due to violation of the inf–sup condition (2.2). Note that this stability issue does not prevent well-posedness of the finite-dimensional discontinuous Galerkin problem. Nor does it prevent convergence to the exact solution as the discretization is refined! The standard analysis technique only requires the weaker condition of *consistency*, that is, that the solution to the original problem satisfies the discrete variational form for each discrete test function.

Through three increasingly complex examples of Friedrichs systems (§ 3–§ 5), the aim of our contribution is to introduce well-posed variational formulations in the sense of (1.6), in which boundary and initial conditions are imposed weakly, similarly as in discontinuous Galerkin methods. The first example, steady advection, constitutes something of a blueprint for the other examples, which is why it is treated in some detail. The analysis of the first example is facilitated by the fact that the solution in this case can be defined directly in the graph space $W$, in contrast to the second example, an elliptic equation written as a first order system, which requires a restricted solution space $V$. However, the analysis of the elliptic problem is simplified by the fact that the graph space can be characterized as a Cartesian product of standard Sobolev spaces. This simplification is not available in the last and most complex example, the acoustic wave equation, which involves inhomogeneous boundary as well as initial conditions. Moreover, in the first two examples, the symmetric part of the operator satisfies, in addition to boundedness condition (1.2), also the following coercivity condition in $L$,

$$c\|\phi\|_L \leq \|(T + \widetilde{T})\phi\|_L \qquad \forall \phi \in \mathscr{D}, \tag{1.7}$$

for some $c > 0$. This condition fails to hold in the last example.

*Remark* 1.3. Although our treatment, as well as our notation, is inspired by EGC [8], we use a slightly different operator formalism, similar to the one used by Wloka [26], for instance. We consistently consider weak differential operators; that is, a derivative is a distributional derivative that can be represented as a locally integrable function, and we typically regard differential operators as bounded linear operators between different spaces instead of as densely-defined unbounded operators. In practice this difference in formalism is nonconsequential, as pointed out by Antonic and Burazin [1, § 2, p. 1697].

The first two examples are similar to examples treated by Ern and Guermond [7] and Ern, Guermond & Caplain [8], except that our theory is variational, which means, for instance, that inhomogeneous boundary conditions are straightforward to handle. The third example below is addressed by Burazin & Erceg [4,

§ 3.3], who treat initial–boundary-value problems using the theory of Ern, Guermond & Caplain [8] for the spatial operator together with the semigroup approach for the time evolution. Our treatment differs in that it is variational, and in that we use a space–time formalism, in which time and space directions are treated on an equal footing. Again, our variational approach makes it straightforward to handle inhomogeneous initial as well as inhomogeneous boundary conditions of characteristic type, also called impedance boundary conditions, which are imposed weakly in the problem statement.

The contributions discussed above [8,11,13,21] aim for a *general theory* of Friedrichs systems, which is *not* the intention here. Rather, we address specific (initial–)boundary-value problems for operators characterized by property (1.2) and employ closely related variational formulations in order to specify precisely in what sense the (initial–)boundary-value problem is set and to establish well-posedness in this sense. We believe that having access to true variational formulations also of Friedrichs-type systems is in itself of interest and closes a "gap" in the classical Hilbert-space approach to the analysis of partial differential equations. Moreover, the variational forms presented below constitute variations of the ones used for discontinuous Galerkin discretizations, and may therefore serve as an inspiration for the development of new numerical methods for Friedrichs systems.

## 2 Well-posedness of variational problems

The well-posedness of variational problem (1.6) is characterized by the following theorem, attributed to Nečas [17].

**Theorem 2.1.** *For real Hilbert spaces $V$ and $L$, let $a$ be a continuous bilinear form on $L \times V$ and $l$ a continuous linear functional on $L$. The variational problem to find $u \in V$ such that $a(v, u) = l(v) \; \forall v \in L$ has a unique solution satisfying*

$$\|u\| \leq \frac{1}{\alpha} \|l\| \tag{2.1}$$

*for some $\alpha > 0$, if and only if the following two conditions hold:*

(i) $\exists \alpha > 0$ *such that, for each $u \in V$,*

$$\sup_{\substack{v \in L \\ v \neq 0}} \frac{a(v, u)}{\|v\|} \geq \alpha \|u\|. \tag{2.2}$$

(ii) *If $v \in L$ satisfies*

$$a(v, u) = 0 \qquad \forall u \in V \tag{2.3}$$

*then $v = 0$.*

The continuous bilinear form defines a bounded linear operator from $V$ to the dual of $L$. Condition (2.2) implies that the operator has a trivial null space and a closed range, and condition (2.3) that it is surjective. Ern and Guermond [6, Thm. 2.6] formulate and prove theorem 2.1 in the more general setting of Banach spaces.

The analysis of variational problem (1.6) is simplified if $L = V$. In particular, a sufficient condition for properties (2.2) and (2.3) is that the bilinear form is strongly coercive. The theorem for this case is known as the Lax–Milgram lemma.

## 3 Example 1: steady advection

A standard model problem for first-order hyperbolic problems is the advection–reaction problem

$$\boldsymbol{\beta} \cdot \nabla u + \rho u = f \qquad \text{in } \Omega, \tag{3.1a}$$
$$u = g \qquad \text{on } \Gamma_-, \tag{3.1b}$$

posed on an open, bounded, and connected Lipschitz domain $\Omega \subset \mathbb{R}^d$, $d = 2$ or $3$. We assume that $\boldsymbol{\beta} \in W^{1,\infty}(\Omega)^d$ and, for simplicity of exposition, that $\nabla \cdot \boldsymbol{\beta} = 0$ and that $\forall \boldsymbol{x} \in \bar{\Omega}$,

$$\rho(\boldsymbol{x}) \geq \rho_0 > 0. \tag{3.2}$$

The boundary $\partial\Omega$ comprises the parts

$$\Gamma_- = \{\, x \in \partial\Omega \mid n \cdot \beta < 0 \,\} \quad \text{(inflow)}, \tag{3.3a}$$
$$\Gamma_+ = \{\, x \in \partial\Omega \mid n \cdot \beta > 0 \,\} \quad \text{(outflow)}, \tag{3.3b}$$
$$\Gamma_0 = \{\, x \in \partial\Omega \mid n \cdot \beta = 0 \,\} \quad \text{(tangential flow)}. \tag{3.3c}$$

One way to generate a variational formulation of problem (3.1) is by a least-squares approach, through which problem (3.1) essentially will be reformulated into an equivalent second-order, anisotropic diffusion problem, as discussed, for instance in Azerad's Ph. D. thesis [2, Ch. 5]. An analogous approach has also been proposed by Bourhrara [3] for the neutron transport equation. However, our aim is to devise variational formulations for the equations in their original Friedrichs-type form. For this, the starting point will be the variational forms that are used within the framework of discontinuous Galerkin methods.

## 3.1 The Discontinuous Galerkin method

The Discontinuous Galerkin (DG) methods for hyperbolic equations was introduced by Reed & Hill [22] and first analyzed by Lesaint & Raviart [16] for model problem (3.1). We will briefly discuss how the method is constructed, since the variational formulation (3.4), from which the DG method can be developed, serves as the starting point also for our approach.

First, let $V_h$ be a finite-dimensional space of weakly differentiable functions—the weak differentiability will later be relaxed—and define the following variational problem.

$$\begin{aligned} &\text{Find } u_h \in V_h \text{ such that} \\ &a_0(v_h, u_h) = l_0(v_h) \qquad \forall v_h \in V_h, \end{aligned} \tag{3.4}$$

where

$$a_0(v_h, u_h) = \int_\Omega v_h (\beta \cdot \nabla u_h + \rho u_h) - \int_{\Gamma_-} n \cdot \beta v_h u_h, \tag{3.5a}$$

$$l_0(v_h) = \int_\Omega v_h f - \int_{\Gamma_-} n \cdot \beta v_h g. \tag{3.5b}$$

*Remark* 3.1. For integrals without "free" variables, like the ones in definitions (3.5), we will in this article not include a measure symbol (such as $dV$ or $dS$), since the type of measure will be clear from the domain of integration.

Note that variational problem (3.4) is *consistent*, that is, $a_0(v_h, u) = l_0(v_h)$, where $u$ is a sufficiently smooth solution (somehow obtained) of boundary-value problem (3.1). Moreover, note that boundary condition (3.1b) is *weakly* imposed, that is, it is not incorporated in the definition of the space but assigned in the variational expression on the same footing as the differential equation in the interior.

The system matrix resulting from problem (3.4) is positive definite since, by choosing $v_h = u_h$ and applying integration by parts,

$$\begin{aligned} a_0(u_h, u_h) &= \int_\Omega u_h (\beta \cdot \nabla u_h + \rho u_h) - \int_{\Gamma_-} n \cdot \beta u_h^2 \\ &= \frac{1}{2} \int_{\partial\Omega} n \cdot \beta u_h^2 + \int_\Omega \rho u_h^2 - \int_{\Gamma_-} n \cdot \beta u_h^2 \\ &= \int_\Omega \rho u_h^2 + \frac{1}{2} \int_{\partial\Omega} |n \cdot \beta| u_h^2 \geq \rho_0 \int_\Omega u_h^2 + \frac{1}{2} \int_{\partial\Omega} |n \cdot \beta| u_h^2 > 0 \qquad \forall u_h \neq 0. \end{aligned} \tag{3.6}$$

It thus follows that system (3.4) is solvable for any data $f \in L^2(\Omega)$, $g \in L^2(\Gamma_-)$.

In spite of the solvability, it turns out that the stability property (3.6) is too weak to obtain accurate approximations. Therefore—and this is the basic feature of DG methods—the continuity of the functions are relaxed to a space of *piecewise* polynomials defined on a triangulation of the domain. Through a bilinear form $a_{\text{DG}}$, inter-element continuity over the edges of the mesh is imposed weakly in the same way as boundary condition (3.1b) is assigned in variational problem (3.4). This procedure leads to the improved stability property

$$a_{\text{DG}}(u_h, u_h) = \int_\Omega \rho u_h^2 + \frac{1}{2} \sum_{K \in \mathscr{T}_h} \int_{\partial K_-} |\boldsymbol{n} \cdot \boldsymbol{\beta}| [\![u_h]\!]^2 + \frac{1}{2} \int_{\Gamma_+} |\boldsymbol{n} \cdot \boldsymbol{\beta}| u_h^2, \qquad (3.7)$$

where $\mathscr{T}_h$ is the set of elements in the mesh, $\partial K_-$ the inflow (cf. definition (3.3a)) portion of the boundary of element $K$, and $[\![u_h]\!]$ the local jump of $u_h$ over the element boundary.

An observation of relevance for what will follow is that the second term in the right-hand side of expression (3.7) can be interpreted as a replacement for the second term in the square of the graph norm

$$\|u\|^2 = \int_\Omega \rho u^2 + \int_\Omega (\boldsymbol{\beta} \cdot \nabla u)^2 \qquad (3.8)$$

associated with operator $\boldsymbol{\beta} \cdot \nabla$.

### 3.2 An ill-posed variational formulation

To set the stage for the later development, it is instructive to see what happens when naively generalizing variational formulation (3.4) to the original, infinite-dimensional boundary-value problem (3.1). Thus, define operator $T = \boldsymbol{\beta} \cdot \nabla + \rho$ and the graph space

$$W = \left\{ u \in L^2(\Omega) \mid Tu \in L^2(\Omega) \right\}, \qquad (3.9)$$

equipped with the norm

$$\|u\|_W = \left( \int_\Omega \left( \rho u^2 + (\boldsymbol{\beta} \cdot \nabla u)^2 \right) \right)^{1/2}, \qquad (3.10)$$

which is equivalent to the graph norm for $T$. As in the discrete case, we choose the space of test functions also as $W$, which leads to the following variational problem.

$$\begin{aligned} &\text{Find } u \in W \text{ such that} \\ &a_0(v, u) = l_0(v) \qquad \forall v \in W. \end{aligned} \qquad (3.11)$$

This variational problem is a particular example of a class of variational formulations discussed by Ern & Guermond [7, Eq. (2.23)], for which they show, in their theorem 2.8, that the solution constructed with their (nonvariational) method will be a unique solution to the variational problem. However, the problem is that variational problem (3.11) in itself is not well posed; the operator defined by $a_0$ does not have a closed range, which means that condition (2.2) will be violated, and the solution will not depend continuously on data. Jensen [13, § 1.9, BVP 2] introduces a similar variational form, with the difference that the space of test functions is a subspace of $H^1(\Omega)$. This formulation suffers from the same shortcoming as problem (3.11).

*Remark* 3.2. In their remark 2.3, Ern & Guermond allude to this problem by stressing that the variational problem does not induce an isomorphism between $W$ and $W'$.

To demonstrate the ill-posedness of problem (3.11), consider the opposite to statement (2.2) applied to $a_0$. That is, for each $\alpha > 0$, there is a $u \in W$ such that

$$\sup_{\substack{v \in W \\ v \neq 0}} \frac{a_0(v, u)}{\|v\|_W} < \alpha \|u\|_W. \qquad (3.12)$$

We will construct such an element $u \in W$. Let $\alpha > 0$ be given and let $(u_n)_{n \in \mathbb{Z}^+}$ be a sequence in $H_0^1(\Omega) \subset W$. Then, using integration by parts and the Cauchy–Schwarz inequality, we find that there is a constant $C > 0$ such that for each element $u_n$ in the sequence and for each $v \in W$,

$$\begin{aligned} a_0(v, u_n) &= \int_\Omega v \boldsymbol{\beta} \cdot \nabla u_n + \int_\Omega \rho v u_n = -\int_\Omega u_n \boldsymbol{\beta} \cdot \nabla v + \int_\Omega \rho u_n v \\ &\leq \|u_n\|_{L^2(\Omega)} \|\boldsymbol{\beta} \cdot \nabla v\|_{L^2(\Omega)} + \|\rho^{1/2} u_n\|_{L^2(\Omega)} \|\rho^{1/2} v\|_{L^2(\Omega)} \\ &\leq C \|u_n\|_{L^2(\Omega)} \|v\|_W, \end{aligned} \qquad (3.13)$$

from which it follows that
$$\sup_{\substack{v \in W \\ v \neq 0}} \frac{a_0(v, u_n)}{\|v\|_W} \leq C \|u_n\|_{L^2(\Omega)}. \tag{3.14}$$

Now choose the sequence to be bounded in $L^2(\Omega)$ but unbounded in $H_0^1(\Omega)$. From inequality (3.14) then follows that there is a $N \in \mathbb{Z}^+$ such that
$$\sup_{\substack{v \in W \\ v \neq 0}} \frac{a_0(v, u_N)}{\|v\|_W} \leq C \|u_N\|_{L^2(\Omega)} \leq C \sup_n \|u_n\|_{L^2(\Omega)} < \alpha \|u_N\|_W, \tag{3.15}$$

where the last two inequalities follows from the boundedness and unboundedness of the sequence in $L^2(\Omega)$ and $H_0^1(\Omega)$, respectively. The element $u = u_N$ thus satisfies inequality (3.12). An example of such a sequence when $\Omega$ is the unit square and $\boldsymbol{\beta} = (1, 0)$ is given by the elements $u_n(x, y) = \sin n\pi x \sin \pi y$.

### 3.3 A well-posed variational formulation

Thus, the bilinear form $a_0$ is not the right choice for a well-posed variational formulation of problem (3.1). To arrive at another formulation, the first crucial observation from inspection of the original problem is that it seems natural to provide the input data $(f, g)$ in the Cartesian product space $L = L^2(\Omega) \times L^2(\Gamma_-; |\boldsymbol{n} \cdot \boldsymbol{\beta}|)$. This choice suggests that we could identify $L' = L$ and use $L$ also as the space of test functions. We thus introduce the *test function tuple* $\hat{v} = (v, v_-) \in L$ and equip $L$ with the norm
$$\|\hat{v}\|_L = \left( \int_\Omega \rho v^2 + \int_{\Gamma_-} |\boldsymbol{n} \cdot \boldsymbol{\beta}| v_-^2 \right)^{1/2} = \left( \|\rho^{1/2} v\|_{L^2(\Omega)}^2 + \|v_-\|_{L^2(\Gamma_-; |\boldsymbol{n} \cdot \boldsymbol{\beta}|)}^2 \right)^{1/2}. \tag{3.16}$$

Moreover, following the framework laid out in § 1, we require the solution space $V$ to be a subspace of the graph space $W$ in definition (3.9). We will assign the boundary conditions weakly through an integral, similarly as in $a_0$. Therefore, due to our choice of $L$, we need to require the solution space $V$ to possess traces in $L^2(\Gamma_-; |\boldsymbol{n} \cdot \boldsymbol{\beta}|)$. In general, continuous trace maps of functions in $W$ can only be defined into $H^{-1/2}(\partial\Omega)$. However, when $\text{dist}(\Gamma_-, \Gamma_+) > 0$, the trace is continuous into $L^2(\partial\Omega; |\boldsymbol{n} \cdot \boldsymbol{\beta}|)$ [7, Lemma 3.1]. Thus, for this example we will assume that $\text{dist}(\Gamma_-, \Gamma_+) > 0$, which means that we can choose the space of solutions simply as $V = W$.

The bilinear form and the linear functional will then become, for $\hat{v} = (v, v_-) \in L, u \in V$,
$$a(\hat{v}, u) = \int_\Omega v T u - \int_{\Gamma_-} \boldsymbol{n} \cdot \boldsymbol{\beta} v_- u, \qquad l(\hat{v}) = \int_\Omega v f - \int_{\Gamma_-} \boldsymbol{n} \cdot \boldsymbol{\beta} v_- g, \tag{3.17}$$

where, as before, $T = \boldsymbol{\beta} \cdot \nabla + \rho$. The variational formulation of boundary-value problem (3.1) can then be stated as follows.
$$\begin{aligned} &\text{Find } u \in V = W \text{ such that} \\ &a(\hat{v}, u) = l(\hat{v}) \qquad \forall \hat{v} \in L. \end{aligned} \tag{3.18}$$

The well-posedness proof of problem (3.18) will also refer to the formal adjoint operator $\widetilde{T}u = -\boldsymbol{\beta} \cdot \nabla u + \rho u$. We note that the graph space associated with $\widetilde{T}$ is also $W$. Complementary to $L$, we also define the space $L^* = L^2(\Omega) \times L^2(\Gamma_+; |\boldsymbol{n} \cdot \boldsymbol{\beta}|)$. Functions $\hat{v} = (v, v_+) \in L^*$ are provided with the norm
$$\|\hat{v}\|_{L^*} = \left( \|\rho^{1/2} v\|_{L^2(\Omega)}^2 + \|v_+\|_{L^2(\Gamma_+; |\boldsymbol{n} \cdot \boldsymbol{\beta}|)}^2 \right)^{1/2}. \tag{3.19}$$

The bilinear form
$$a^*(\hat{v}, u) = \int_\Omega v \widetilde{T} u + \int_{\Gamma_+} \boldsymbol{n} \cdot \boldsymbol{\beta} v_+ u \tag{3.20}$$
is continuous on $L^* \times V$, and is the adjoint of $a$ in the sense that, for each $u, v \in V$,
$$a^*(\hat{u}, v) = a(\hat{v}, u), \tag{3.21}$$
where $\hat{u} = (u, \text{tr}_+ u)$ and $\hat{v} = (v, \text{tr}_- v)$, and where $\text{tr}_+$ and $\text{tr}_-$ are the trace maps into $L^2(\Gamma_-; |\boldsymbol{n} \cdot \boldsymbol{\beta}|)$ and $L^2(\Gamma_+; |\boldsymbol{n} \cdot \boldsymbol{\beta}|)$, respectively, that is,
$$\text{tr}_- \in \mathscr{L}(V, L^2(\Gamma_-; |\boldsymbol{n} \cdot \boldsymbol{\beta}|)), \qquad \text{tr}_+ \in \mathscr{L}(V, L^2(\Gamma_+; |\boldsymbol{n} \cdot \boldsymbol{\beta}|)). \tag{3.22}$$

Now, we are prepared to establish the following well-posedness result.

**Theorem 3.3.** *With $a$ and $l$ as in definition (3.17), assuming bound (3.2) and that* $\mathrm{dist}(\Gamma_-, \Gamma_+) > 0$, *the variational problem to find* $u \in V = W$ *such that*

$$a(\hat{v}, u) = l(\hat{v}) \qquad \forall \hat{v} \in L \tag{3.23}$$

*has a unique solution satisfying*

$$\|u\| \leq 2\|l\|. \tag{3.24}$$

*Proof.* Since $\mathrm{dist}(\Gamma_-, \Gamma_+) > 0$, the trace maps (3.22) are well defined [7, lemma 3.1]. From the Cauchy–Schwarz inequality follows then that the bilinear forms $a$, $a^*$ and the linear functional $l$ are continuous. Thus, in order to apply theorem 2.1, we need to show that conditions (2.2) and (2.3) are satisfied.

*Condition (2.2):* The condition is satisfied for $u = 0$. For $u \in V \setminus \{0\}$, choose $\hat{v} = \hat{u} = (u, \mathrm{tr}_- u) \in L$. Integrating by parts, we then find

$$\begin{aligned}a(\hat{u}, u) &= \int_\Omega u[(\boldsymbol{\beta} \cdot \nabla)u + \rho u] - \int_{\Gamma_-} \boldsymbol{n} \cdot \boldsymbol{\beta} u^2 = \frac{1}{2}\int_{\partial\Omega} |\boldsymbol{n} \cdot \boldsymbol{\beta}|u^2 + \frac{1}{2}\int_\Omega \rho u^2 + \frac{1}{2}\int_\Omega \rho u^2 \\ &\geq \|\hat{u}\|_L \|\rho^{1/2} u\|_{L^2(\Omega)},\end{aligned} \tag{3.25}$$

from which follows that

$$\sup_{\hat{v} \in L \setminus \{0\}} \frac{a(\hat{v}, u)}{\|\hat{v}\|_L} \geq \|\rho^{1/2} u\|_{L^2(\Omega)}. \tag{3.26}$$

An analogous calculation (now choosing $\hat{v} = (u, \mathrm{tr}_+ u) \in L^*$) reveals the same bound for $a^*$,

$$\sup_{\hat{v} \in L^* \setminus \{0\}} \frac{a^*(\hat{v}, u)}{\|\hat{v}\|_{L^*}} \geq \|\rho^{1/2} u\|_{L^2(\Omega)}. \tag{3.27}$$

Moreover, for each $u \in V$ such that $Tu \neq 0$, choosing $\hat{w} = (Tu, 0) \in L$, we find that

$$\sup_{\hat{v} \in L \setminus \{0\}} \frac{a(\hat{v}, u)}{\|\hat{v}\|_L} \geq \frac{a(\hat{w}, u)}{\|\hat{w}\|_L} = \|Tu\|_{L^2(\Omega)}. \tag{3.28}$$

From bounds (3.26) and (3.28) it follows that for each $u \in V \setminus \{0\}$,

$$2 \sup_{\hat{v} \in L \setminus \{0\}} \frac{a(\hat{v}, u)}{\|\hat{v}\|_L} \geq \|\rho^{1/2} u\|_{L^2(\Omega)} + \|Tu\|_{L^2(\Omega)} \geq \|u\|_V, \tag{3.29}$$

We have thus verified condition (2.2) in theorem 2.1 with $\alpha = 1/2$.

*Condition (2.3):* Let $\hat{v} = (v, v_-) \in L$ such that

$$a(\hat{v}, u) = \int_\Omega v(\boldsymbol{\beta} \cdot \nabla u + \rho u) - \int_{\Gamma_-} \boldsymbol{n} \cdot \boldsymbol{\beta} v_- u = 0 \qquad \forall u \in V. \tag{3.30}$$

Choosing $u = (\phi, 0)$ in equation (3.30), where $\phi \in C_0^1(\Omega)$, we find that

$$\int_\Omega v T\phi = 0 \qquad \forall \phi \in C_0^1(\Omega). \tag{3.31}$$

That is, by the definition of weak derivative,

$$\widetilde{T} v = 0, \tag{3.32}$$

so that, in particular, $v \in W$, which means that $v$ admits $L^2(\partial\Omega)$ traces. We may therefore, for any $u \in C^1(\bar{\Omega})$, integrate expression (3.30) by parts to obtain

$$\begin{aligned}0 &= \int_\Omega v(\boldsymbol{\beta} \cdot \nabla u + \rho u) - \int_{\Gamma_-} \boldsymbol{n} \cdot \boldsymbol{\beta} v_- u = \int_{\partial\Omega} \boldsymbol{n} \cdot \boldsymbol{\beta} vu + \int_\Omega u \underbrace{\widetilde{T} v}_{=0} - \int_{\Gamma_-} \boldsymbol{n} \cdot \boldsymbol{\beta} v_- u \\ &= \int_{\Gamma_+} |\boldsymbol{n} \cdot \boldsymbol{\beta}| uv + \int_{\Gamma_-} |\boldsymbol{n} \cdot \boldsymbol{\beta}| u(v_- - v) = 0.\end{aligned} \tag{3.33}$$

Choosing $u$ as functions whose trace on $\partial\Omega$ has compact support in $\Gamma_+$ and $\Gamma_-$, respectively, and by density, we conclude that

$$\text{tr}_+ v = 0, \tag{3.34a}$$
$$\text{tr}_- v = v_-. \tag{3.34b}$$

Expressions (3.32) and (3.34a) substituted in definition (3.20) yield that $a^*(w, v) = 0 \ \forall w \in L^*$, which means that $v = 0$ by inequality (3.27). Expression (3.34b) implies then that also $v_- = 0$, which finally shows that $\hat{v} = (v, v_-) = 0$, which verifies also condition (2.3) and thereby, by theorem 2.1, shows well-posedness of variational problem (3.18). □

## 4 Example 2: an elliptic model problem

Our second example concerns the following boundary-value problem for a vector field $\boldsymbol{u}$ and a scalar field $p$,

$$\boldsymbol{u} + \nabla p = \boldsymbol{f}_1 \quad \text{in } \Omega, \tag{4.1a}$$
$$p + \nabla \cdot \boldsymbol{u} = f_2 \quad \text{in } \Omega, \tag{4.1b}$$
$$\frac{1}{2}(1-\alpha)p - \frac{1}{2}(1+\alpha)\boldsymbol{n} \cdot \boldsymbol{u} = g \quad \text{on } \partial\Omega, \tag{4.1c}$$

which constitutes a first-order-system formulation of the scalar second-order elliptic problem

$$-\Delta p + p = f \quad \text{in } \Omega, \tag{4.2a}$$
$$\frac{1}{2}(1-\alpha)p + \frac{1}{2}(1+\alpha)\frac{\partial p}{\partial n} = g \quad \text{on } \partial\Omega. \tag{4.2b}$$

We assume the domain $\Omega$ to be open, bounded, connected, and Lipschitz. Moreover, the function $\alpha \in L^\infty(\partial\Omega)$ is assumed to satisfy, for some $\alpha_M \in [0, 1)$,

$$\text{ess. im } \alpha \subset [-\alpha_M, \alpha_M]. \tag{4.3}$$

That is, the interpolation in the Robin-type boundary condition (4.1c) is not allowed anywhere to reduce to a pure Dirichlet ($\alpha = -1$) or Neumann ($\alpha = 1$) condition on $p$.

*Remark* 4.1. The reason for the restriction in $\alpha$ is that the inf–sup constant of the variational formulation will in our proof turn out to be proportional to $1 - \alpha_M$.

*Remark* 4.2. Homogeneous pure Dirichlet and Neumann conditions can be handled, due to the characterization in lemma 4.3 of the graph space, by incorporating these strongly in the components of the solution vector. We choose to ignore this case for simplicity of exposition.

Equations (4.1a) and (4.1b) can be written in the block operator form

$$\mathbf{T}\boldsymbol{\xi} = \mathbf{f} \tag{4.4}$$

where

$$\boldsymbol{\xi} = \begin{bmatrix} \boldsymbol{\xi}_1 \\ \xi_2 \end{bmatrix}, \qquad \mathbf{T} = \begin{bmatrix} \boldsymbol{I} & \nabla \\ \nabla \cdot & 1 \end{bmatrix}, \qquad \mathbf{f} = \begin{bmatrix} \boldsymbol{f}_1 \\ f_2 \end{bmatrix}. \tag{4.5}$$

Note the blocking of the rows in $\boldsymbol{\xi}$ and $\mathbf{f}$ in a vector $\boldsymbol{\xi}_1$ (in the sense of a first-order tensor of dimension $d = 2$ or 3, the space dimension) and a scalar $\xi_2$. Consequently, the first column of matrix $\mathbf{T}$ contains operators acting on vector fields and the second column operators that act on scalar fields.

Proceeding similarly as for the first example, we introduce the graph space associated with block operator $\mathbf{T}$ from which a solution space will be extracted,

$$W = \left\{ \boldsymbol{\xi} \in L^2(\Omega)^{d+1} \mid \mathbf{T}\boldsymbol{\xi} \in L^2(\Omega)^{d+1} \right\}. \tag{4.6}$$

However, for this particular example, there is a more elementary characterization of $W$.

**Lemma 4.3.** *It holds that*

$$W = \left\{ \boldsymbol{\xi} = [\boldsymbol{\xi}_1, \xi_2]^T \in L^2(\Omega)^{d+1} \mid \boldsymbol{\xi}_1 \in H(\mathrm{div}; \Omega), \xi_2 \in H^1(\Omega) \right\} \tag{4.7}$$

*Proof.* Definition (4.7) corresponds to the graph norm of $\mathbf{A} = \begin{bmatrix} 0 & \nabla \\ \nabla \cdot & 0 \end{bmatrix}$. Since $\mathbf{T} = \mathbf{I} + \mathbf{A}$, the conclusion follows from that the graph norm of $\mathbf{A}$ and $\mathbf{I} + \mathbf{A}$ are equivalent. □

Characterization (4.7) enables the integration-by-parts formula

$$\begin{aligned}
\int_\Omega \boldsymbol{\eta}^T \mathbf{T} \boldsymbol{\xi} &= \int_\Omega \left[ \boldsymbol{\eta}_1 \cdot (\boldsymbol{\xi}_1 + \nabla \xi_2) + \eta_2 (\xi_2 + \nabla \cdot \boldsymbol{\xi}_1) \right] \\
&= \langle \gamma_{\boldsymbol{n}} \boldsymbol{\eta}_1, \xi_2 \rangle + \langle \gamma_{\boldsymbol{n}} \boldsymbol{\xi}_1, \eta_2 \rangle + \int_\Omega \boldsymbol{\xi}^T \tilde{\mathbf{T}} \boldsymbol{\eta} \qquad \forall \boldsymbol{\eta}, \boldsymbol{\xi} \in W,
\end{aligned} \tag{4.8}$$

where

$$\tilde{\mathbf{T}} = \begin{pmatrix} \mathbf{I} & -\nabla \\ -\nabla \cdot & 1 \end{pmatrix} \tag{4.9}$$

is the formal adjoint of $\mathbf{T}$, where $\langle \cdot, \cdot \rangle$ denotes the duality pairing on $H^{-1/2}(\partial\Omega) \times H^{1/2}(\partial\Omega)$, and where $\gamma_{\boldsymbol{n}} \in \mathscr{L}\big(H(\mathrm{div}, \Omega), H^{-1/2}(\partial\Omega)\big)$ is the continuous extension of the trace operator that for $\boldsymbol{u} \in C^1(\overline{\Omega})^d$ satisfies $\gamma_{\boldsymbol{n}} \boldsymbol{u} = \boldsymbol{n} \cdot \boldsymbol{u}|_{\partial\Omega}$. In particular, for $\boldsymbol{\eta} = \boldsymbol{\xi}$, formula (4.8) reduces to

$$\int_\Omega \boldsymbol{\xi}^T \mathbf{T} \boldsymbol{\xi} = \int_\Omega |\boldsymbol{\xi}|^2 + \langle \gamma_{\boldsymbol{n}} \boldsymbol{\xi}_1, \xi_2 \rangle, \qquad \forall \boldsymbol{\xi} \in W. \tag{4.10}$$

*Remark* 4.4. The graph space corresponding to operator $\tilde{\mathbf{T}}$ is identical to $W$. That is, in addition to definition (4.6), it holds that

$$W = \left\{ \boldsymbol{\xi} \in L^2(\Omega)^{d+1} \mid \tilde{\mathbf{T}} \boldsymbol{\xi} \in L^2(\Omega)^{d+1} \right\}. \tag{4.11}$$

In order to generalize the approach of § 3.3 to system (4.1), we first notice that it seems reasonable to provide data to system (4.1) as a tuple $(\mathbf{f}, g)$ of interior data $\mathbf{f} = [\boldsymbol{f}_1, f_2]^T \in L^2(\Omega)^{d+1}$ and boundary data $g \in L^2(\partial\Omega)$. Consequently, we therefore define the space of test functions as the Cartesian product space

$$L = L^2(\Omega)^{d+1} \times L^2(\partial\Omega). \tag{4.12}$$

For elements in $L$, we will use the same tuple notation as for the data, that is, $\hat{\boldsymbol{\eta}} = (\boldsymbol{\eta}, \eta_R)$, where $\boldsymbol{\eta} = [\boldsymbol{\eta}_1, \eta_2]^T \in L^2(\Omega)^{d+1}$ and $\eta_R \in L^2(\partial\Omega)$, and provide the norm

$$\|\boldsymbol{\eta}\|_L = \left( \int_\Omega \left( |\boldsymbol{\eta}_1|^2 + \eta_2^2 \right) + \int_{\partial\Omega} |\eta_R|^2 \right)^{1/2}. \tag{4.13}$$

Note that the first element $\boldsymbol{\eta}$ in the test-function tuple will have the same block structure as the elements in $W$ and correspond to the interior data vector $\mathbf{f}$. The second element $\eta_R$ in the test-function tuple correspond to the scalar boundary data $g$. Associated with the boundary condition, we introduce the trace map

$$\mathrm{tr}_\alpha \boldsymbol{\xi} = \frac{1}{\sqrt{2}} \left[ (1-\alpha)\xi_2 - (1+\alpha)\boldsymbol{n} \cdot \boldsymbol{\xi}_1 \right]\big|_{\partial\Omega}, \tag{4.14}$$

defined for $\boldsymbol{\xi} \in C^1(\overline{\Omega})^{d+1}$. For $\hat{\boldsymbol{\eta}} = (\boldsymbol{\eta}, \eta_R) \in L$ and $\boldsymbol{\xi} \in V$, where $V \subset W$ is a suitable solution space, below defined so that $\mathrm{tr}_\alpha$ can be continuously extended to $\mathscr{L}(V, L^2(\partial\Omega))$, we define

$$a(\hat{\boldsymbol{\eta}}, \boldsymbol{\xi}) = \int_\Omega \boldsymbol{\eta}^T \mathbf{T} \boldsymbol{\xi} + \int_{\partial\Omega} \eta_R \, \mathrm{tr}_\alpha \boldsymbol{\xi}, \tag{4.15a}$$

$$l(\hat{\boldsymbol{\eta}}) = \int_\Omega \boldsymbol{\eta}^T \mathbf{f} + \sqrt{2} \int_{\partial\Omega} \eta_R g. \tag{4.15b}$$

The issue is now to define a suitable space of solutions $V \subset W$ to render a well-posed variational formulation. By characterization (4.7), we see that a restriction of $\boldsymbol{\xi}_1$ is needed to admit traces $\mathrm{tr}_\alpha$ in $L^2(\partial\Omega)$, since $H(\text{div};\Omega)$-functions generally admits normal traces only in $H^{-1/2}(\partial\Omega)$. Therefore, we introduce the following strict subspace of $H(\text{div};\Omega)$

$$U = \left\{ \boldsymbol{u} \in H(\text{div};\Omega) \mid \gamma_{\boldsymbol{n}} \boldsymbol{u} \in L^2(\partial\Omega) \right\}, \tag{4.16}$$

equipped with the inner product

$$(\boldsymbol{u},\boldsymbol{v})_U = \int_\Omega (\boldsymbol{u}\cdot\boldsymbol{v} + \nabla\cdot\boldsymbol{u}\nabla\cdot\boldsymbol{v}) + \int_\Gamma \gamma_{\boldsymbol{n}}\boldsymbol{u}\,\gamma_{\boldsymbol{n}}\boldsymbol{v}. \tag{4.17}$$

To show that $U$, as well as the solution space of example 3 in § 5, is a Hilbert space, we will rely on the following general result.

**Theorem 4.5.** *Let $X$, $Y$, and $Z$ be Banach spaces such that $Y \subset Z$ with continuous embedding, and let $A : X \to Z$ be a bounded linear operator. Then the space*

$$X_Y = \{\, x \in X \mid Ax \in Y \,\}, \tag{4.18}$$

*with norm*

$$\|x\|_{X_Y} = \bigl(\|x\|_X^2 + \|Ax\|_Y^2\bigr)^{1/2}, \tag{4.19}$$

*is a Banach space continuously embedded in $X$.*

*Proof.* Since $\|x\|_X \leq \|x\|_{X_Y}$, $X_Y$ embeds continuously into $X$. It remains to show that $X_Y$ is complete. Let the sequence $(x_n)_{n\in\mathbb{Z}^+}$ be Cauchy in $X_Y$. By the continuous embedding, $(x_n)_{n\in\mathbb{Z}^+}$ is Cauchy also in $X$, so there is an $x^*$ such that $x_n \to x^*$ in $X$. Moreover, $(Ax_n)_{n\in\mathbb{Z}^+}$ is Cauchy in $Y$, so there is a $y^* \in Y$ such that $Ax_n \to y^*$ in $Y$. Thus, $X_Y$ will be complete if $Ax^* = y$. By the continuous embedding $Y \subset Z$, $Ax_n \to y^*$ also in $Z$. Since also $Ax_n \to Ax^*$ in $Z$, by continuity of $A$, uniqueness of limits yields that $Ax^* = y$ and $X_Y$ is thus complete. □

Choosing $X = H(\text{div};\Omega)$, $Y = L^2(\Omega)$, $Z = H^{-1/2}(\partial\Omega)$, $X_Y = U$, and $A = \gamma_{\boldsymbol{n}}$, theorem 4.5 implies following result.

**Lemma 4.6.** *The space $U$ is a Hilbert space continuously embedded in $H(\text{div};\Omega)$.*

Now we are ready to define the solution space as the Hilbert space

$$V = \left\{\, \boldsymbol{\xi} = (\xi_1, \xi_2) \in W \mid \xi_1 \in U \,\right\}, \tag{4.20}$$

equipped with norm

$$\|\boldsymbol{\xi}\|_V = \left(\|\boldsymbol{\xi}\|_W^2 + \|\gamma_{\boldsymbol{n}}\xi_1\|_{L^2(\partial\Omega)}^2\right)^{1/2}. \tag{4.21}$$

Since $\gamma_{\boldsymbol{n}}$ maps functions in $U$ into $L^2(\partial\Omega)$, integration-by-parts formula (4.8) can, in the particular case of $\boldsymbol{\eta} \in W$, $\boldsymbol{\xi} \in V$ be simplified and written

$$\int_\Omega \boldsymbol{\eta}^T \mathbf{T}\boldsymbol{\xi} = \langle \gamma_{\boldsymbol{n}}\eta_1, \xi_2 \rangle + \int_{\partial\Omega} \gamma_{\boldsymbol{n}}\xi_1 \eta_2 + \int_\Omega \boldsymbol{\xi}^T \widetilde{\mathbf{T}}\boldsymbol{\eta}, \tag{4.22}$$

and, in particular, for $\boldsymbol{\xi} \in V$,

$$\int_\Omega \boldsymbol{\xi}^T \mathbf{T}\boldsymbol{\xi} = \int_\Omega |\boldsymbol{\xi}|^2 + \int_{\partial\Omega} \gamma_{\boldsymbol{n}}\xi_1 \xi_2. \tag{4.23}$$

The variational problem corresponding to boundary-value problem (4.1) can then be formulated in standard form.

$$\begin{aligned}&\text{Find } \boldsymbol{\xi} \in V \text{ such that} \\ &a(\hat{\boldsymbol{\eta}}, \boldsymbol{\xi}) = l(\hat{\boldsymbol{\eta}}) \qquad \forall \hat{\boldsymbol{\eta}} \in L.\end{aligned} \tag{4.24}$$

We will also utilize the bilinear form

$$a^*(\hat{\boldsymbol{\eta}}, \boldsymbol{\xi}) = \int_\Omega \boldsymbol{\eta}^T \widetilde{\mathbf{T}}\boldsymbol{\xi} + \int_{\partial\Omega} \eta_R \, \text{tr}_\alpha^* \boldsymbol{\xi}, \tag{4.25}$$

where the trace operator
$$\mathrm{tr}_\alpha^* \boldsymbol{\xi} = \frac{1}{\sqrt{2}}\big[(1-\alpha)\xi_2 + (1+\alpha)\boldsymbol{n}\cdot\boldsymbol{\xi}_1\big]\big|_{\partial\Omega} \tag{4.26}$$

also extends continuously to a bounded operator from $V$ into $L^2(\partial\Omega)$. The form $a^*$ is adjoint to $a$ in the sense that for $\boldsymbol{\xi} \in V$ and $\boldsymbol{\eta} \in V$,
$$a(\hat{\boldsymbol{\eta}}, \boldsymbol{\xi}) = a^*(\hat{\boldsymbol{\xi}}, \boldsymbol{\eta}) \tag{4.27}$$

where
$$\hat{\boldsymbol{\eta}} = (\boldsymbol{\eta}, \mathrm{tr}_0\,\boldsymbol{\eta}), \qquad \hat{\boldsymbol{\xi}} = (\boldsymbol{\xi}, \mathrm{tr}_0^*\,\boldsymbol{\xi}). \tag{4.28}$$

To prove the inf-sup condition (2.2), we first establish the following bounds.

**Lemma 4.7.** *For each $\boldsymbol{\xi} \in V$, it holds that*
$$\sup_{\hat{\boldsymbol{\eta}} \in L\setminus\{0\}} \frac{a(\hat{\boldsymbol{\eta}}, \boldsymbol{\xi})}{\|\hat{\boldsymbol{\eta}}\|} \geq \frac{1-\alpha_M}{2}\left(\|\boldsymbol{\xi}\|^2_{L^2(\Omega)^{d+1}} + \|\gamma_{\boldsymbol{n}}\boldsymbol{\xi}_1\|^2_{L^2(\partial\Omega)}\right)^{1/2}, \tag{4.29a}$$
$$\sup_{\hat{\boldsymbol{\eta}} \in L\setminus\{0\}} \frac{a^*(\hat{\boldsymbol{\eta}}, \boldsymbol{\xi})}{\|\hat{\boldsymbol{\eta}}\|} \geq \frac{1-\alpha_M}{2}\left(\|\boldsymbol{\xi}\|^2_{L^2(\Omega)^{d+1}} + \|\gamma_{\boldsymbol{n}}\boldsymbol{\xi}_1\|^2_{L^2(\partial\Omega)}\right)^{1/2}. \tag{4.29b}$$

*Proof.* Due to the bilinearity of $a$ and $a^*$, the conditions hold for $\boldsymbol{\xi} = \boldsymbol{0}$. Thus, let $\boldsymbol{\xi} \in V$ be nonzero and define
$$\hat{\boldsymbol{\xi}} = (\boldsymbol{\xi}, \mathrm{tr}_0\,\boldsymbol{\xi}) \in L. \tag{4.30}$$

The conclusion (4.29a) then follows from the calculation

$$\begin{aligned}
a(\hat{\boldsymbol{\xi}}, \boldsymbol{\xi}) &= \int_\Omega \boldsymbol{\xi}^T \mathbf{T}\boldsymbol{\xi} + \frac{1}{2}\int_{\partial\Omega}\Big[(\xi_2 - \gamma_{\boldsymbol{n}}\boldsymbol{\xi}_1)((1-\alpha)\xi_2 - (1+\alpha)\gamma_{\boldsymbol{n}}\boldsymbol{\xi}_1)\Big] \\
&= \int_\Omega |\boldsymbol{\xi}|^2 + \int_{\partial\Omega}\xi_2\,\gamma_{\boldsymbol{n}}\boldsymbol{\xi}_1 + \frac{1}{2}\int_{\partial\Omega}\Big[(1-\alpha)\xi_2^2 + (1+\alpha)(\gamma_{\boldsymbol{n}}\boldsymbol{\xi}_1)^2\Big] - \int_{\partial\Omega}\xi_2\,\gamma_{\boldsymbol{n}}\boldsymbol{\xi}_1 \\
&\geq \int_\Omega |\boldsymbol{\xi}|^2 + \frac{1-\alpha_M}{2}\int_{\partial\Omega}(\xi_2^2 + (\gamma_{\boldsymbol{n}}\boldsymbol{\xi}_1)^2) \\
&= \left[\int_\Omega |\boldsymbol{\xi}|^2 + \frac{1-\alpha_M}{2}\int_{\partial\Omega}(\xi_2^2 + (\gamma_{\boldsymbol{n}}\boldsymbol{\xi}_1)^2)\right]^{\frac{1}{2}}\left[\int_\Omega |\boldsymbol{\xi}|^2 + \frac{1-\alpha_M}{2}\int_{\partial\Omega}(\xi_2^2 + (\gamma_{\boldsymbol{n}}\boldsymbol{\xi}_1)^2)\right]^{\frac{1}{2}} \\
&\geq \left[\int_\Omega |\boldsymbol{\xi}|^2 + \frac{1-\alpha_M}{2}\int_{\partial\Omega}(\gamma_{\boldsymbol{n}}\boldsymbol{\xi}_1)^2\right]^{\frac{1}{2}}\left[\int_\Omega |\boldsymbol{\xi}|^2 + \frac{1-\alpha_M}{2}\int_{\partial\Omega}\frac{1}{2}(\xi_2 - \gamma_{\boldsymbol{n}}\boldsymbol{\xi}_1)^2\right]^{\frac{1}{2}} \\
&\geq \frac{1-\alpha_M}{2}\left(\|\boldsymbol{\xi}\|^2_{L^2(\Omega)^{d+1}} + \|\gamma_{\boldsymbol{n}}\boldsymbol{\xi}_1\|^2_{L^2(\partial\Omega)}\right)^{1/2}\|\hat{\boldsymbol{\xi}}\|_L,
\end{aligned} \tag{4.31}$$

where the second equality follows from integration-by-parts formula (4.23), the first inequality from the bound (4.3) on $\alpha$, and where in the second inequality, we have neglected $\xi_2^2$ in the first factor and used
$$a^2 + b^2 = \frac{1}{2}\big((a-b)^2 + (a+b)^2\big) \geq \frac{1}{2}(a-b)^2 \tag{4.32}$$

for $a = \xi_2$ and $b = \gamma_{\boldsymbol{n}}\boldsymbol{\xi}_1$ in the second factor. The dual conclusion (4.29b) follows by an analogous calculation on $a^*$ using test function
$$\hat{\boldsymbol{\xi}} = (\boldsymbol{\xi}, \mathrm{tr}_0^*\,\boldsymbol{\xi}) \in L. \tag{4.33}$$
□

With the help of lemma 4.7, the required inf–sup condition is straightforward to show.

**Lemma 4.8.** *For each $\boldsymbol{\xi} \in V$, it holds that*
$$\sup_{\hat{\boldsymbol{\eta}} \in L\setminus\{0\}} \frac{a(\hat{\boldsymbol{\eta}}, \boldsymbol{\xi})}{\|\hat{\boldsymbol{\eta}}\|} \geq \frac{1-\alpha_M}{3}\|\boldsymbol{\xi}\|_V. \tag{4.34}$$

*Proof.* Let $\boldsymbol{\xi} \in V$. If $\mathbf{T}\boldsymbol{\xi} = \mathbf{0}$, then the conclusion follows from lemma 4.7. We may thus assume that $\mathbf{T}\boldsymbol{\xi} \neq \mathbf{0}$. Define $\widetilde{\boldsymbol{\eta}} = (\mathbf{T}\boldsymbol{\xi}, 0) \in L$. Then

$$\|\widetilde{\boldsymbol{\eta}}\|_L = \|\mathbf{T}\boldsymbol{\xi}\|_{L^2(\Omega)^{d+1}}, \tag{4.35}$$

and

$$a(\widetilde{\boldsymbol{\eta}}, \boldsymbol{\xi}) = \int_\Omega (\mathbf{T}\boldsymbol{\xi})^T \mathbf{T}\boldsymbol{\xi} = \|\mathbf{T}\boldsymbol{\xi}\|^2_{L^2(\Omega)^{d+1}}, \tag{4.36}$$

which implies that

$$\sup_{\widehat{\boldsymbol{\eta}} \in L \setminus \{\mathbf{0}\}} \frac{a(\widehat{\boldsymbol{\eta}}, \boldsymbol{\xi})}{\|\widehat{\boldsymbol{\eta}}\|} \geq \frac{a(\widetilde{\boldsymbol{\eta}}, \boldsymbol{\xi})}{\|\widetilde{\boldsymbol{\eta}}\|} = \|\mathbf{T}\boldsymbol{\xi}\|_{L^2(\Omega)^{d+1}}. \tag{4.37}$$

Inequalities (4.29a) and (4.37) imply that

$$\left(\frac{2}{1-\alpha_M} + 1\right) \sup_{\widehat{\boldsymbol{\eta}} \in L \setminus \{\mathbf{0}\}} \frac{a(\widehat{\boldsymbol{\eta}}, \boldsymbol{\xi})}{\|\widehat{\boldsymbol{\eta}}\|} \geq \left(\|\boldsymbol{\xi}\|^2_{L^2(\Omega)^{d+1}} + \|\gamma_n \boldsymbol{\xi}_1\|^2_{L^2(\partial\Omega)}\right)^{1/2} + \|\mathbf{T}\boldsymbol{\xi}\|_{L^2(\Omega)^{d+1}} \tag{4.38}$$

$$\geq \|\boldsymbol{\xi}\|_V,$$

from which the conclusion follows. $\square$

What is now left is to show surjectivity.

**Lemma 4.9.** *If $\widehat{\boldsymbol{\eta}} \in L$ such that*

$$a(\widehat{\boldsymbol{\eta}}, \boldsymbol{\xi}) = 0 \qquad \forall \boldsymbol{\xi} \in V, \tag{4.39}$$

*then $\widehat{\boldsymbol{\eta}} = \mathbf{0}$.*

*Proof.* By definition (4.15a), condition (4.39) reads $\widehat{\boldsymbol{\eta}} = (\boldsymbol{\eta}, \eta_R) \in L$, where $\boldsymbol{\eta} = [\boldsymbol{\eta}_1, \eta_2]$, such that

$$\int_\Omega \boldsymbol{\eta}^T \mathbf{T}\boldsymbol{\xi} + \int_{\partial\Omega} \eta_R \operatorname{tr}_\alpha \boldsymbol{\xi} = 0 \qquad \forall \boldsymbol{\xi} = [\boldsymbol{\xi}_1, \xi_2]^T \in V, \tag{4.40}$$

from which it follows that

$$\int_\Omega \boldsymbol{\eta}^T \mathbf{T}\boldsymbol{\phi} = 0 \qquad \forall \boldsymbol{\phi} \in C_0^1(\Omega)^{d+1}; \tag{4.41}$$

that is,

$$\widetilde{\mathbf{T}}\boldsymbol{\eta} = \mathbf{0} \tag{4.42}$$

by the definition of weak derivative. We conclude thus that, trivially, $\boldsymbol{\eta} \in W$ (Remark 4.4), which means that we may integrate the first term in equation (4.40) by parts, using formula (4.22), to obtain

$$\langle \gamma_n \boldsymbol{\eta}_1, \xi_2 \rangle_{\partial\Omega} + \int_{\partial\Omega} \gamma_n \boldsymbol{\xi}_1 \, \eta_2 + \frac{1}{\sqrt{2}} \int_{\partial\Omega} \eta_R \big((1-\alpha)\xi_2 - (1+\alpha)\gamma_n \boldsymbol{\xi}_1\big) = 0 \quad \forall [\boldsymbol{\xi}_1, \xi_2]^T \in V. \tag{4.43}$$

In particular, for $\boldsymbol{\xi}_1 = \mathbf{0}$, $\xi_2 \in H^1(\Omega)$, we find that

$$\left\langle \gamma_n \boldsymbol{\eta}_1 + \frac{1-\alpha}{\sqrt{2}} \eta_R, \gamma \xi_2 \right\rangle = 0 \qquad \forall \xi_2 \in H^1(\Omega), \tag{4.44}$$

where $\gamma$ is the trace map of $H^1(\Omega)$ onto $H^{1/2}(\partial\Omega)$. Since $\gamma$ is surjective, it follows that

$$\gamma_n \boldsymbol{\eta}_1 + \frac{1-\alpha}{\sqrt{2}} \eta_R = 0, \tag{4.45}$$

and, in particular, that $\gamma_n \boldsymbol{\eta}_1 \in L^2(\partial\Omega)$ (since $\eta_R$ is in $L^2(\partial\Omega)$).

Choosing $\xi_2 = 0$, $\boldsymbol{\xi}_1 \in C^\infty(\overline{\Omega})^d$ in equation (4.43), we find that

$$\int_{\partial\Omega} \gamma_n \boldsymbol{\xi}_1 \left(\eta_2 - \frac{1+\alpha}{\sqrt{2}} \eta_R\right) = 0 \qquad \forall \boldsymbol{\xi}_1 \in C^\infty(\overline{\Omega})^d, \tag{4.46}$$

from which it follows that
$$\eta_2 - \frac{1+\alpha}{\sqrt{2}} \eta_R = 0. \tag{4.47}$$

Multiplying expressions (4.45) and (4.47) with $(1+\alpha)$ and $(1-\alpha)$, respectively, and adding, it follows that
$$\frac{1}{\sqrt{2}}\bigl[(1+\alpha)\gamma_{\boldsymbol{n}}\boldsymbol{\eta}_1 + (1-\alpha)\eta_2\bigr] = \operatorname{tr}^*_\alpha \boldsymbol{\eta} = 0. \tag{4.48}$$

Due to expressions (4.42) and (4.48), and by definition (4.25), we see that $a^*(\hat{\boldsymbol{\theta}}, \boldsymbol{\eta}) = 0$ for each $\hat{\boldsymbol{\theta}} \in L$. Property (4.29b) then implies that $\boldsymbol{\eta} = \boldsymbol{0}$. Finally, equation (4.47) yields that also $\eta_R$ vanishes; hence $\hat{\boldsymbol{\eta}} = \boldsymbol{0}$. □

With these results, well-posedness of the variational problem is straightforward to show.

**Theorem 4.10.** *With $a$ and $l$ as in definitions (4.15), where function $\alpha$ satisfies bound (4.3), with space $V$ as in definition (4.20), and with $L = L^2(\Omega)^{d+1} \times L^2(\partial\Omega)$, the variational problem to find $\boldsymbol{\xi} \in V$ such that*
$$a(\hat{\boldsymbol{\eta}}, \boldsymbol{\xi}) = l(\hat{\boldsymbol{\eta}}) \qquad \forall \hat{\boldsymbol{\eta}} \in L \tag{4.49}$$

*has a unique solution satisfying*
$$\|\boldsymbol{\xi}\| \leq \frac{3}{1-\alpha_M} \|l\|. \tag{4.50}$$

*Proof.* By the Cauchy–Schwarz inequality, $a$ and $l$ are continuous on $L \times V$ and $L$, respectively. Well-posedness then follows from theorem 2.1 together with lemmas 4.8 and 4.9. □

## 5  Example 3: the acoustic wave equation

Here we consider the equations of linear acoustics in a still, ideal gas under isentropic conditions. The two previous examples were idealized model problem, templates for basic hyperbolic and elliptic equations, respectively, without the inclusion of appropriate dimensional coefficients that would occur in applications. In contrast, the equations and the spaces discussed here will be presented in a form that respects relevant physical units.

The gaseous medium is characterized by its static density $\rho_0$ and speed of sound $c_0$. In simple situations, the static density and the speed of sound are constant, but in the presence of temperature gradients in the gas, the density and speed of sound will vary spatially. However, typically the medium's static pressure as well as the quantity $\rho_0 c_0^2$, can be regarded as constant, also in the presence of temperature gradients. The constancy of these quantities follows from the linearization of the Euler equations of gas dynamics in the case when body forces acting on the system can be neglected. Motivations for these assumptions and more details on the modeling are given by Rienstra & Hirschberg [23, § 2.4].

The boundary-value problem under consideration here will be

$$\rho_0 \frac{\partial \boldsymbol{u}}{\partial t} + \nabla p = \boldsymbol{f}_1 \qquad \text{in } Q = \Omega \times (0, \tau), \tag{5.1a}$$

$$\frac{\partial p}{\partial t} + \rho_0 c_0^2 \nabla \cdot \boldsymbol{u} = f_2 \qquad \text{in } Q = \Omega \times (0, \tau), \tag{5.1b}$$

$$\frac{1}{2}(p - \rho_0 c_0 \boldsymbol{n} \cdot \boldsymbol{u}) - \frac{\alpha}{2}(p + \rho_0 c_0 \boldsymbol{n} \cdot \boldsymbol{u}) = g \qquad \text{on } \Sigma = \partial\Omega \times (0, \tau), \tag{5.1c}$$

$$\boldsymbol{u} = \boldsymbol{u}_s \qquad p = p_s \qquad \text{on } Q_0 = \Omega \times \{0\}. \tag{5.1d}$$

The unknown quantities are the acoustic velocity and pressure fields $\boldsymbol{u}$ and $p$, and data to the system is provided through the right-hand forcing in equations (5.1a), (5.1b), boundary condition (5.1c), and initial conditions (5.1d). The system (5.1a), (5.1b) constitutes a first-order-systems formulation of the scalar wave equation

$$\frac{\partial^2 p}{\partial t^2} - \nabla \cdot c_0^2 \nabla p = f \qquad \text{in } \Omega \times (0, \tau). \tag{5.2}$$

We assume the domain $\Omega$ to be open, bounded, and connected with a smooth boundary $\partial\Omega$, and that the domain locally is located on one side of its boundary.

*Remark* 5.1. The smoothness property that is assumed of $\partial\Omega$ in the analysis below is that it is $C^1$, due to theorem 5.10, with a normal field $\boldsymbol{n}$ that is Lipschitz, due to the application of the Kirszbraun theorem.

As can be noted already in the formulation (5.1), we adopt a "space–time" formalism: the function spaces will be defined on the space–time cylinder $Q$, a Lipschitz domain whose boundary $\partial Q$ is naturally partitioned as

$$\partial Q = \bar{\Sigma} \cup \bar{Q}_0 \cup \bar{Q}_\tau. \tag{5.3}$$

Moreover, the interpolation function $\alpha \in L^\infty(\Sigma)$ in boundary condition (5.1c) is assumed, analogously as in § 4, to satisfy, for some $\alpha_M \in [0, 1)$,

$$\text{ess. im}\, \alpha \subset [-\alpha_M, \alpha_M]. \tag{5.4}$$

*Remark* 5.2. Note that restriction (5.4) means that boundary condition (5.1c) cannot reduce to a pure Dirichlet condition on $p$ or $\boldsymbol{n} \cdot \boldsymbol{u}$.

Exploiting that $\rho_0 c_0^2$ is constant and that media properties $\rho_0$ and $c_0$ have no time dependency, equations (5.1a) and (5.1b) can be rewritten in the block operator form

$$\mathbf{T}\boldsymbol{\xi} = \mathbf{f}, \tag{5.5}$$

in which

$$\boldsymbol{\xi} = \begin{bmatrix} \boldsymbol{\xi}_1 \\ \xi_2 \end{bmatrix} = \begin{bmatrix} \rho_0 c_0 \boldsymbol{u} \\ p \end{bmatrix}, \qquad \mathbf{T} = \begin{bmatrix} \partial_t \boldsymbol{I} & c_0 \nabla \\ \nabla c_0 \cdot & \partial_t \end{bmatrix}, \qquad \mathbf{f} = \begin{bmatrix} c_0 \boldsymbol{f}_1 \\ f_2 \end{bmatrix}, \tag{5.6}$$

with the same blocking of the $d + 1$ rows of $\boldsymbol{\xi}$ and $\mathbf{f}$ as in § 4. The Cartesian components of operator $\mathbf{T}$ in $d = 3$ are

$$[\mathbf{T}] = \begin{pmatrix} \frac{\partial}{\partial t} & 0 & 0 & c_0 \frac{\partial}{\partial x_1} \\ 0 & \frac{\partial}{\partial t} & 0 & c_0 \frac{\partial}{\partial x_2} \\ 0 & 0 & \frac{\partial}{\partial t} & c_0 \frac{\partial}{\partial x_3} \\ \frac{\partial}{\partial x_1} c_0 & \frac{\partial}{\partial x_2} c_0 & \frac{\partial}{\partial x_3} c_0 & \frac{\partial}{\partial t} \end{pmatrix}. \tag{5.7}$$

The formal adjoint of $\mathbf{T}$ is $\widetilde{\mathbf{T}} = -\mathbf{T}$.

Under the assumptions discussed above, equation (5.5) holds also for a spatially varying speed of sound $c_0$, generated by temperature gradients in the medium. However, from now on, in order to simplify the analysis, we will assume that $c_0$ is constant and positive. The graph space associated with block operator $\mathbf{T}$ will be

$$W = \left\{ \boldsymbol{\xi} \in L^2(Q)^{d+1} \mid \mathbf{T}\boldsymbol{\xi} \in L^2(Q)^{d+1} \right\}, \tag{5.8}$$

equipped with norm

$$\|\boldsymbol{\xi}\|_W = \left( \|\boldsymbol{\xi}\|_{L^2(Q)^{d+1}}^2 + \tau^2 \|\mathbf{T}\boldsymbol{\xi}\|_{L^2(Q)^{d+1}}^2 \right)^{1/2}, \tag{5.9}$$

and we note that the graph space associated with $\widetilde{\mathbf{T}}$ is also $W$. Note that we scale the velocity unknowns so that all components of $\boldsymbol{\xi}$ will have the same dimension (pressure). Moreover, by the inclusion of constants $c_0$ and $\tau$ in the definition of $\mathbf{T}$ and the norm on $W$, all terms that are summed will possess the same dimension. Consequently, in this section, it will be convenient also to equip $H^1(Q)$ with the dimensionally consistent norm

$$\|u\|_{H^1(Q)} = \left[ \|u\|_{L^2(Q)}^2 + \tau^2 \left( \|\partial_t u\|_{L^2(Q)}^2 + \|c_0 \nabla u\|_{L^2(Q)}^2 \right) \right]^{1/2}. \tag{5.10}$$

By the inequality

$$|\boldsymbol{\xi}|^2 + \tau^2 |\mathbf{T}\boldsymbol{\xi}|^2 = |\boldsymbol{\xi}|^2 + \tau^2 |\partial_t \boldsymbol{\xi}_1 + c_0 \nabla \xi_2|^2 + \tau^2 |\partial_t \xi_2 + \nabla \cdot c_0 \boldsymbol{\xi}_1|^2 \\ \leq |\boldsymbol{\xi}|^2 + 2\tau^2 \left[ |\partial_t \boldsymbol{\xi}_1|^2 + |c_0 \nabla \xi_2|^2 + |\partial_t \xi_2|^2 + |c_0 \nabla \boldsymbol{\xi}_1|^2 \right], \tag{5.11}$$

we conclude that

$$\|\boldsymbol{\xi}\|_W \leq \sqrt{2} \|\boldsymbol{\xi}\|_{H^1(Q)^{d+1}}, \tag{5.12}$$

and thus that $H^1(Q)^{d+1} \subset W$. However, in this case, as opposed to the elliptic case of § 4, there is no characterization of $W$ as simple as in lemma 4.3.

To introduce a variational formulation of system (5.1), we will proceed similarly as in previous sections and consider data to the problem as being given in the form of the tuple

$$(\mathbf{f}, g, \boldsymbol{\xi}_s) \in L^2(Q)^{d+1} \times L^2(\Sigma) \times L^2(Q_0)^{d+1}, \tag{5.13}$$

where $\mathbf{f} = [\boldsymbol{f}_1, f_2]^T$ and $\boldsymbol{\xi}_s = [\boldsymbol{u}_s, p_s]^T$. This form of the data suggests a space of test functions of the same form,

$$\hat{\boldsymbol{\eta}} = (\boldsymbol{\eta}, \eta_\Sigma, \boldsymbol{\eta}_s) \in L^2(Q)^{d+1} \times L^2(\Sigma) \times L^2(Q_0)^{d+1} = L, \tag{5.14}$$

which we equip with the norm

$$\|\hat{\boldsymbol{\eta}}\|_L = \left(\|\boldsymbol{\eta}\|^2_{L^2(Q)^{d+1}} + \tau\|\eta_\Sigma\|^2_{L^2(\Sigma)} + \tau\|\boldsymbol{\eta}_s\|^2_{L^2(Q_0)}\right)^{1/2}. \tag{5.15}$$

The components of the test function tuple will be used to enforce the equation system, the boundary conditions, and the initial conditions, respectively.

Associated with the boundary and initial conditions, we introduce the trace maps

$$\mathrm{tr}^\pm_\Sigma \boldsymbol{\xi} = \sqrt{\frac{c_0}{2}}(\xi_2 \pm \boldsymbol{n} \cdot \boldsymbol{\xi}_1)|_\Sigma, \quad \mathrm{tr}_{Q_0} \boldsymbol{\xi} = \boldsymbol{\xi}|_{Q_0}, \quad \mathrm{tr}_{Q_\tau} \boldsymbol{\xi} = \boldsymbol{\xi}|_{Q_\tau}, \tag{5.16}$$

defined for $\boldsymbol{\xi} \in C^1(\overline{Q})^{d+1}$. For $\hat{\boldsymbol{\eta}} \in L$ and $\boldsymbol{\xi} \in V$, where below we will define $V \subset W$ so that the ranges of the trace operators (5.16) continuously extend into $L^2$ spaces, we define

$$a(\hat{\boldsymbol{\eta}}, \boldsymbol{\xi}) = \int_Q \boldsymbol{\eta}^T \mathbf{T}\boldsymbol{\xi} + \int_\Sigma \eta_\Sigma \left(\mathrm{tr}^-_\Sigma \boldsymbol{\xi} - \alpha \, \mathrm{tr}^+_\Sigma \boldsymbol{\xi}\right) + \int_{Q_0} \boldsymbol{\eta}_s^T \left(\mathrm{tr}_{Q_0} \boldsymbol{\xi}\right), \tag{5.17a}$$

$$l(\hat{\boldsymbol{\eta}}) = \int_Q \boldsymbol{\eta}^T \boldsymbol{f} + \sqrt{2c_0} \int_\Sigma \eta_\Sigma g + \int_{Q_0} \boldsymbol{\eta}_s^T \boldsymbol{\xi}_s. \tag{5.17b}$$

We also define, complementary to the space $L$, the space $L^* = L^2(Q^{d+1}) \times L^2(\Sigma) \times L^2(Q_\tau)$, equipped with norm

$$\|\hat{\boldsymbol{\eta}}\|_{L^*} = \left(\|\boldsymbol{\eta}\|^2_{L^2(Q)^{d+1}} + \tau\|\eta_\Sigma\|^2_{L^2(\Sigma)} + \tau\|\boldsymbol{\eta}_s\|^2_{L^2(Q_\tau)}\right)^{1/2}, \tag{5.18}$$

and the adjoint bilinear form, to be used in the surjectivity proof,

$$a^*(\hat{\boldsymbol{\eta}}, \boldsymbol{\xi}) = -\int_Q \boldsymbol{\eta}^T \mathbf{T}\boldsymbol{\xi} + \int_\Sigma \eta_\Sigma \left(\mathrm{tr}^+_\Sigma \boldsymbol{\xi} - \alpha \, \mathrm{tr}^-_\Sigma \boldsymbol{\xi}\right) + \int_{Q_\tau} \boldsymbol{\eta}_s^T \left(\mathrm{tr}_{Q_\tau} \boldsymbol{\xi}\right), \tag{5.19}$$

which satisfies $a^*(\hat{\boldsymbol{\xi}}, \boldsymbol{\eta}) = a(\hat{\boldsymbol{\eta}}, \boldsymbol{\xi})$ for $\boldsymbol{\eta}, \boldsymbol{\xi} \in C^1(\overline{Q})^{d+1}$ and

$$\hat{\boldsymbol{\eta}} = (\boldsymbol{\eta}, \mathrm{tr}^-_\Sigma \boldsymbol{\eta}, \mathrm{tr}_{Q_0} \boldsymbol{\eta}), \qquad \hat{\boldsymbol{\xi}} = (\boldsymbol{\xi}, \mathrm{tr}^+_\Sigma \boldsymbol{\xi}, \mathrm{tr}_{Q_\tau} \boldsymbol{\xi}). \tag{5.20}$$

The basic integration-by-parts formula for operator $\mathbf{T}$, repeatedly used in the following, is as follows. Let $\boldsymbol{\eta}, \boldsymbol{\xi} \in C^1(\overline{Q})^{d+1}$. Then

$$\begin{aligned}
\int_Q \boldsymbol{\eta}^T \mathbf{T}\boldsymbol{\xi} &= \int_Q [\eta_1, \eta_2] \begin{bmatrix} \partial_t \boldsymbol{\xi}_1 + c_0 \nabla \xi_2 \\ \partial_t \xi_2 + \nabla \cdot c_0 \boldsymbol{\xi}_1 \end{bmatrix} \\
&= \int_Q \left(\boldsymbol{\eta}_1 \cdot \partial_t \boldsymbol{\xi}_1 + c_0 \boldsymbol{\eta}_1 \cdot \nabla \xi_2 + \eta_2 \partial_t \xi_2 + c_0 \eta_2 \nabla \cdot \boldsymbol{\xi}_1\right) \\
&= \int_{Q_\tau} \boldsymbol{\eta}^T \boldsymbol{\xi} - \int_{Q_0} \boldsymbol{\eta}^T \boldsymbol{\xi} + \int_\Sigma c_0(\boldsymbol{n} \cdot \boldsymbol{\eta}_1 \xi_2 + \eta_2 \boldsymbol{n} \cdot \boldsymbol{\xi}_1) + \int_Q \boldsymbol{\xi}^T \underbrace{\widetilde{\mathbf{T}} \boldsymbol{\eta}}_{=-\mathbf{T}\boldsymbol{\eta}} \\
&= \int_{\partial Q} \boldsymbol{\eta}^T \mathbf{T}_\nu \boldsymbol{\xi} - \int_Q \boldsymbol{\xi}^T \mathbf{T}\boldsymbol{\eta},
\end{aligned} \tag{5.21}$$

where
$$\mathbf{T}_\nu \boldsymbol{\xi} = \begin{cases} c_0 [\boldsymbol{n}\xi_2, \boldsymbol{n}\cdot\boldsymbol{\xi}_1]^T & \text{on } \Sigma, \\ \boldsymbol{\xi} & \text{on } Q_\tau, \\ -\boldsymbol{\xi} & \text{on } Q_0. \end{cases} \tag{5.22}$$

As proven by Jensen [13, Thm. 4], for instance, the following density property holds for domains (like $Q$) that possess the segment property.

**Theorem 5.3.** *The space $C^\infty(\overline{Q})^{d+1}$ is dense in $W$.*

Due to this property, we show next that the basic integration-by-parts formula (5.21) extends to $\boldsymbol{\eta} \in H^1(Q)^{d+1}$ and $\boldsymbol{\xi} \in W$.

**Lemma 5.4.** *Assume that $Q$ is a space-time cylinder satisfying the segment property. The trace map $\mathbf{T}_\nu$, defined in expression (5.22) for functions $\boldsymbol{\xi} \in C^1(\overline{Q})^{d+1}$, extends continuously to $\mathscr{L}(W, H^{-1/2}(\partial Q)^{d+1})$, and integration-by-parts formula (5.21) extends to $\boldsymbol{\eta} \in H^1(Q)^{d+1}$ and $\boldsymbol{\xi} \in W$, so that*

$$\int_Q \boldsymbol{\eta}^T \mathbf{T}\boldsymbol{\xi} = \langle \mathbf{T}_\nu \boldsymbol{\xi}, \gamma\boldsymbol{\eta}\rangle_{H^{1/2}(\partial Q)^{d+1}} - \int_Q \boldsymbol{\xi}^T \mathbf{T}\boldsymbol{\eta}, \tag{5.23}$$

*where $\gamma \in \mathscr{L}(H^1(Q), H^{1/2}(\partial Q))$ denotes the boundary trace map. Moreover, the bound*

$$\langle \mathbf{T}_\nu \boldsymbol{\xi}, \boldsymbol{\psi}\rangle_{H^{1/2}(\partial Q)^{d+1}} \le \frac{4}{\tau} \|\gamma_*\| \|\boldsymbol{\psi}\|_{H^{1/2}(\partial Q)^{d+1}} \|\boldsymbol{\xi}\|_W \tag{5.24}$$

*holds for $\boldsymbol{\xi} \in W$ and $\boldsymbol{\psi} \in H^{1/2}(\partial Q)^{d+1}$, and where $\gamma_* \in \mathscr{L}(H^{1/2}(\partial Q)^{d+1}, H^1(Q)^{d+1})$ denotes a right inverse of $\gamma$.*

*Remark* 5.5. Analogously as in § 4, $\langle \cdot, \cdot\rangle_{H^{1/2}(\partial Q)}$ denotes the duality pairing on $H^{-1/2}(\partial Q) \times H^{1/2}(\partial Q)$.

*Proof.* Let $\boldsymbol{\xi} \in C^1(\overline{Q})^{d+1}$. By density of $C^1(\overline{Q})^{d+1}$ in $H^1(Q)^{d+1}$, continuity of the $H^1(Q)^{d+1}$ boundary trace, and estimate (5.12), we find that integration-by-parts formula (5.21) extends to $\boldsymbol{\eta} \in H^1(Q)^{d+1}$, so that

$$\int_Q \boldsymbol{\eta}^T \mathbf{T}\boldsymbol{\xi} = \int_{\partial Q} (\gamma\boldsymbol{\eta})^T \mathbf{T}_\nu \boldsymbol{\xi} - \int_Q \boldsymbol{\xi}^T \mathbf{T}\boldsymbol{\eta}, \tag{5.25}$$

from which it follows that

$$\int_{\partial Q} (\gamma\boldsymbol{\eta})^T \mathbf{T}_\nu \boldsymbol{\xi} = \int_Q \boldsymbol{\eta}^T \mathbf{T}\boldsymbol{\xi} + \int_Q \boldsymbol{\xi}^T \mathbf{T}\boldsymbol{\eta} \le \frac{2}{\tau} \|\boldsymbol{\eta}\|_W \|\boldsymbol{\xi}\|_W \le \frac{4}{\tau} \|\boldsymbol{\eta}\|_{H^1(Q)^{d+1}} \|\boldsymbol{\xi}\|_W, \tag{5.26}$$

using inequality (5.12) in the last step.

Let $\boldsymbol{\psi} \in H^{1/2}(\partial Q)^{d+1}$ and $\gamma_*$ a continuous right inverse of $\gamma$. Then inequality (5.26) gives that

$$\begin{aligned}\int_{\partial Q} \boldsymbol{\psi}^T \mathbf{T}_\nu \boldsymbol{\xi} &= \int_Q (\gamma_*\boldsymbol{\psi})^T \mathbf{T}\boldsymbol{\xi} + \int_Q \boldsymbol{\xi}^T \mathbf{T}(\gamma_*\boldsymbol{\psi}) \\ &\le \frac{4}{\tau} \|\gamma_*\boldsymbol{\psi}\|_{H^1(Q)^{d+1}} \|\boldsymbol{\xi}\|_W \le \frac{4}{\tau} \|\gamma_*\| \|\boldsymbol{\psi}\|_{H^{1/2}(\partial Q)^{d+1}} \|\boldsymbol{\xi}\|_W.\end{aligned} \tag{5.27}$$

By inequality (5.27), we find that

$$\|\mathbf{T}_\nu \boldsymbol{\xi}\|_{H^{-1/2}(\partial Q)^{d+1}} = \sup_{\substack{\boldsymbol{\psi} \in H^{1/2}(\partial Q)^{d+1} \\ \boldsymbol{\psi} \ne 0}} \frac{\int_{\partial Q} \boldsymbol{\psi}^T \mathbf{T}_\nu \boldsymbol{\xi}}{\|\boldsymbol{\psi}\|_{H^{1/2}(\partial Q)^{d+1}}} \le \frac{4}{\tau} \|\gamma_*\| \|\boldsymbol{\xi}\|_W. \tag{5.28}$$

By density theorem 5.3, it follows that $\mathbf{T}_\nu$ extends continuously to $\mathscr{L}(W, H^{-1/2}(\partial Q)^{d+1})$. Integration-by-parts formula (5.25) therefore extends to formula (5.23) and bound (5.27) to bound (5.24). □

The following technical lemma will be used in the trace theorems below. We suspect this lemma or variants thereof to be known. However, we have failed to find a suitable reference and therefore provide a proof in the appendix.

**Lemma 5.6.** *Let $\boldsymbol{h} \in C^{0,\mu}(\partial Q)^n$ be a Hölder continuous function with exponent $\mu \in (1/2, 1]$. Then there is a constant $C$ such that for any $u \in H^{1/2}(\partial Q)$,*

$$\|\boldsymbol{h} u\|_{H^{1/2}(\partial Q)^n} \leq C \|u\|_{H^{1/2}(\partial Q)}. \tag{5.29}$$

More precisely, lemma 5.6 will be applied with $\boldsymbol{h} = \mathbf{h}_{\pm} \in \text{Lip}(\partial Q)^{d+1}$ defined by

$$\mathbf{h}_{\pm} = \frac{1}{\sqrt{2c_0}} \begin{bmatrix} \boldsymbol{n}_* \\ \pm 1 \end{bmatrix}, \tag{5.30}$$

where $\boldsymbol{n}_*$ denotes the Lipschitz extension to $\partial Q$ of the normal field $\boldsymbol{n}$ on $\Sigma$, which exists due to the Kirszbraun theorem [14, Thm. 5.2.2]. Note that the multiplier $\mathbf{h}_{\pm}$ is defined so that, for $\boldsymbol{\xi} \in C^1(\bar{Q})^{d+1}$, $\mathbf{h}_{\pm}^T \mathbf{T}_\nu \boldsymbol{\xi} = \text{tr}_\Sigma^\pm \boldsymbol{\xi}$ on $\Sigma$ (recall definitions (5.22) and (5.16)).

Now we are ready to prove that trace maps (5.16) are well defined also for arguments in $W$. Each of these traces turns out to be definable with ranges in the dual of the so-called Lions–Magenes space [24, § 33] $H_{00}^{1/2}(\Gamma)$, where $\Gamma$ is $Q_0$, $Q_\tau$, or $\Sigma$. This space can be defined as follows. For $g \in L^2(\Gamma)$, denote by $g_*$ the extension by zero of $g$ to all of $\partial Q$. Then

$$H_{00}^{1/2}(\Gamma) = \left\{ g \in L^2(\Gamma) \mid \exists u \in H^1(Q) \text{ such that } \gamma u = g_* \right\}, \tag{5.31}$$

provided with the quotient norm

$$\|g\|_{H_{00}^{1/2}(\Gamma)} = \inf_{\substack{u \in H^1(Q) \\ \gamma v = g_*}} \|u\|_{H^1(Q)}. \tag{5.32}$$

It holds that $H_{00}^{1/2}(\Gamma) \subset H^{1/2}(\Gamma)$, with continuous embedding. The typical use case for this space is, like here, to characterize boundary conditions on $\Gamma$ as residing in the dual space $H_{00}^{1/2}(\Gamma)'$.

**Lemma 5.7.** *The trace maps $\text{tr}_{Q_0}$, $\text{tr}_{Q_\tau}$, and $\text{tr}_\Sigma^\pm$, defined in expression (5.16) for $\boldsymbol{\xi} \in C^1(\bar{Q})^{d+1}$, extend continuously to $\mathscr{L}(W, (H_{00}^{1/2}(Q_0)^{d+1})')$, $\mathscr{L}(W, (H_{00}^{1/2}(Q_\tau)^{d+1})')$, and $\mathscr{L}(W, H_{00}^{1/2}(\Sigma)')$, respectively. Moreover,*

*(i) for $\boldsymbol{\xi} \in W$ and $\boldsymbol{\psi} \in H_{00}^{1/2}(Q_0)^{d+1}$,*

$$\langle \text{tr}_{Q_0} \boldsymbol{\xi}, \boldsymbol{\psi} \rangle_{H_{00}^{1/2}(Q_0)^{d+1}} = \langle \mathbf{T}_\nu \boldsymbol{\xi}, -\boldsymbol{\psi}_* \rangle_{H^{1/2}(\partial Q)^{d+1}}, \tag{5.33}$$

*where, $\boldsymbol{\psi}_* \in H^{1/2}(\partial Q)^{d+1}$ denotes the extension by zero of $\boldsymbol{\psi}$, $\gamma$ the trace map of $H^1(Q)$ onto $H^{1/2}(\partial Q)$ and $\gamma_*$ a continuous right inverse of $\gamma$.*

*(ii) for $\boldsymbol{\xi} \in W$ and $\boldsymbol{\psi} \in H_{00}^{1/2}(Q_\tau)^{d+1}$,*

$$\langle \text{tr}_{Q_\tau} \boldsymbol{\xi}, \boldsymbol{\psi} \rangle_{H_{00}^{1/2}(Q_\tau)^{d+1}} = \langle \mathbf{T}_\nu \boldsymbol{\xi}, \boldsymbol{\psi}_* \rangle_{H^{1/2}(\partial Q)^{d+1}}, \tag{5.34}$$

*where $\boldsymbol{\psi}_* \in H^{1/2}(\partial Q)^{d+1}$ denotes the extension by zero of $\boldsymbol{\psi}$.*

*(iii) for $\boldsymbol{\xi} \in W$ and $\psi \in H_{00}^{1/2}(\Sigma)$,*

$$\langle \text{tr}_\Sigma^\pm \boldsymbol{\xi}, \psi \rangle_{H_{00}^{1/2}(\Sigma)} = \langle \mathbf{T}_\nu \boldsymbol{\xi}, \mathbf{h}_\pm \psi_* \rangle_{H^{1/2}(\partial Q)^{d+1}}, \tag{5.35}$$

*where $\psi_* \in H^{1/2}(\partial Q)$ denotes the extension by zero of $\psi$, and $\mathbf{h}_\pm \in \text{Lip}(\partial Q)^{d+1}$ is defined in expression (5.30).*

*Proof.* (i): Let $\boldsymbol{\psi} \in H_{00}^{1/2}(Q_0)^{d+1}$, and let $\boldsymbol{\psi}_* \in H^{1/2}(\partial Q)^{d+1}$ be its extension by zero. Then there exist positive constants $C_1 = 4\|\gamma_*\|/\tau$ and $C_2$ such that, for $\boldsymbol{\xi} \in C^1(\bar{Q})^{d+1}$,

$$\int_{Q_0} \boldsymbol{\psi}^T \operatorname{tr}_{Q_0} \boldsymbol{\xi} = -\int_{\partial Q} \boldsymbol{\psi}_*^T \mathbf{T}_\nu \boldsymbol{\xi} \le C_1 \|\boldsymbol{\psi}_*\|_{H^{1/2}(\partial Q)^{d+1}} \|\boldsymbol{\xi}\|_W = C_1 \|\boldsymbol{\psi}\|_{H^{1/2}(Q_0)^{d+1}} \|\boldsymbol{\xi}\|_W \\ \le C_2 \|\boldsymbol{\psi}\|_{H_{00}^{1/2}(Q_0)^{d+1}} \|\boldsymbol{\xi}\|_W, \qquad (5.36)$$

where definition (5.22) is used in the first equality, the bound (5.24) in the first inequality, and the continuous embedding $H_{00}^{1/2}(Q_0)^{d+1} \subset H^{1/2}(Q_0)^{d+1}$ in the last. Inequality (5.36) implies that

$$\|\operatorname{tr}_{Q_0} \boldsymbol{\xi}\|_{(H_{00}^{1/2}(Q_0)^{d+1})'} = \sup_{\substack{\boldsymbol{\psi} \in H_{00}^{1/2}(Q_0)^{d+1} \\ \boldsymbol{\psi} \ne 0}} \frac{\int_{Q_0} \boldsymbol{\psi}^T \operatorname{tr}_{Q_0} \boldsymbol{\xi}}{\|\boldsymbol{\psi}\|_{H_{00}^{1/2}(Q_0)^{d+1}}} \le C_2 \|\boldsymbol{\xi}\|_W, \qquad (5.37)$$

from which it follows that $\operatorname{tr}_{Q_0}$ extends to $\mathcal{L}(W, (H_{00}^{1/2}(Q_0)^{d+1})')$ by density theorem 5.3, and thus that the first equality in expression (5.36) extends to identity (5.33).

(ii): The conclusions follows by analogous arguments as in (i) by considering $\boldsymbol{\psi} \in H_{00}^{1/2}(Q_\tau)^{d+1}$.

(iii): Let $\psi \in H_{00}^{1/2}(\Sigma)$ and let $\psi_* \in H^{1/2}(\partial Q)$ be its extension by zero. By lemma 5.6, $\mathbf{h}_\pm \psi_* \in H^{1/2}(\partial Q)^{d+1}$ and there exists a constant $C_h$ such that $\|\mathbf{h}_\pm \psi_*\|_{H^{1/2}(\partial Q)^{d+1}} \le C_h \|\psi_*\|_{H^{1/2}(\partial Q)}$, where $\mathbf{h}_\pm \in \operatorname{Lip}(\partial Q)^{d+1}$ is defined by expression (5.30) so that, for $\boldsymbol{\xi} \in C^1(\bar{Q})^{d+1}$, $\mathbf{h}_\pm^T \mathbf{T}_\nu \boldsymbol{\xi} = \operatorname{tr}_\Sigma^\pm \boldsymbol{\xi}$ on $\Sigma$. Then there are positive constants $C_1 = 4\|\gamma_*\|/\tau$, $C_2$, and $C_3$ such that, for $\boldsymbol{\xi} \in C^1(\bar{Q})^{d+1}$,

$$\int_\Sigma \psi \operatorname{tr}_\Sigma^\pm \boldsymbol{\xi} = \int_{\partial Q} (\mathbf{h}_\pm \psi_*)^T \mathbf{T}_\nu \boldsymbol{\xi} \le C_1 \|\mathbf{h}_\pm \psi_*\|_{H^{1/2}(\partial Q)^{d+1}} \|\boldsymbol{\xi}\|_W \\ \le C_2 \|\psi_*\|_{H^{1/2}(\partial Q)} \|\boldsymbol{\xi}\|_W = C_2 \|\psi\|_{H^{1/2}(\Sigma)} \|\boldsymbol{\xi}\|_W \le C_3 \|\psi\|_{H_{00}^{1/2}(\Sigma)} \|\boldsymbol{\xi}\|_W, \qquad (5.38)$$

where the bound (5.24) is used in the first inequality, lemma (5.6) in the second, and, in the last inequality, the continuous embedding $H_{00}^{1/2}(\Sigma) \subset H^{1/2}(\Sigma)$. Inequality (5.38), together with the definition of the dual norm and density theorem 5.3 implies that $\operatorname{tr}_\Sigma^\pm$ extends to $\mathcal{L}(W, H_{00}^{1/2}(\Sigma)')$, which in turn implies that the first equality in expression (5.38) extends to identity (5.35). □

Since $L^2(Q_0) \subset (H_{00}^{1/2}(Q_0))'$ and $L^2(\Sigma) \subset (H_{00}^{1/2}(\Sigma))'$, lemma 5.7 implies that the space

$$V = \left\{ \boldsymbol{\xi} \in W \mid \operatorname{tr}_{Q_0} \boldsymbol{\xi} \in L^2(Q_0)^{d+1} \text{ and } \operatorname{tr}_\Sigma^- \boldsymbol{\xi} \in L^2(\Sigma) \right\}, \qquad (5.39)$$

equipped with the norm

$$\|\boldsymbol{\xi}\|_V = \left( \|\boldsymbol{\xi}\|_W^2 + \tau \|\operatorname{tr}_\Sigma^- \boldsymbol{\xi}\|_{L^2(\Sigma)}^2 + \tau \|\operatorname{tr}_{Q_0} \boldsymbol{\xi}\|_{L^2(Q_0)}^2 \right)^{1/2}, \qquad (5.40)$$

is well defined. Note that the norm also can be written

$$\|\boldsymbol{\xi}\|_V = \left( \|(\boldsymbol{\xi}, \operatorname{tr}_\Sigma^- \boldsymbol{\xi}, \operatorname{tr}_{Q_0} \boldsymbol{\xi})\|_L^2 + \tau^2 \|\mathbf{T}\boldsymbol{\xi}\|_{L^2(Q)^{d+1}}^2 \right)^{1/2}, \qquad (5.41)$$

a form that will be utilized in the proof of the inf–sup condition of $a$. We will see that $V$ is a suitable solution space. First we establish that $V$ is a Hilbert space.

**Lemma 5.8.** *The space $V$ is a Hilbert space continuously embedded in $W$.*

*Proof.* Since the norm (5.40) can be derived from a inner product, it remains just to show that $V$ is complete. In order to use theorem 4.5, we identify $X = W$, $Y = L^2(\Sigma) \times L^2(Q_0)$, $Z = (H_{00}^{1/2}(\Sigma)^{d+1})' \times H_{00}^{1/2}(Q_0)'$, $X_Y = V$, and $A = \operatorname{tr}_\Sigma^- \times \operatorname{tr}_{Q_0}$. Lemma 5.7 shows the required properties of $A$, and the fact that $H_{00}^{1/2}(\Gamma) \subset H^{1/2}(\Gamma) \subset L^2(\Gamma)$, for $\Gamma = \Sigma$ or $Q_0$, with continuous embeddings, implies by duality the continuous embedding $Y \subset Z$. Completeness then follows from theorem 4.5. □

The proof of the "extended" trace property in lemma 5.11 below uses density theorem 5.10, whose proof in turn relies on the following density theorem due to Rauch [21, Theorem 8].

**Theorem 5.9.**

(i) The space $C^1(\bar{Q})^{d+1} \cap V_0$ is a dense subspace of $V_0 = \{\boldsymbol{\xi} \in W \mid \operatorname{tr}_{Q_0} \boldsymbol{\xi} = 0 \text{ and } \operatorname{tr}_\Sigma^- \boldsymbol{\xi} = 0\}$.

(ii) The space $C^1(\bar{Q})^{d+1} \cap \widetilde{V}_0$ is a dense subspace of $\widetilde{V}_0 = \{\boldsymbol{\xi} \in W \mid \operatorname{tr}_{Q_\tau} \boldsymbol{\xi} = 0 \text{ and } \operatorname{tr}_\Sigma^+ \boldsymbol{\xi} = 0\}$.

A prerequisite for theorem 5.9 is that $\Sigma$ is a characteristic surface of constant multiplicity; that is, the dimension of $\{\boldsymbol{\eta} \in \mathbb{R}^{d+1} \mid \mathbf{T}_\nu(x)\boldsymbol{\eta} = 0\}$ is independent of $x$ on $\Sigma$, which is true for $\mathbf{T}_\nu$ in expression (5.22).

**Theorem 5.10.** *The space $C^1(\bar{Q})^{d+1}$ is a dense subspace of $V$.*

*Proof.* By lemma 5.8, $V$ is a Hilbert space. The inner product generating the norm (5.40) is, for $\boldsymbol{\eta}, \boldsymbol{\xi} \in V$,

$$(\boldsymbol{\eta}, \boldsymbol{\xi})_V = (\boldsymbol{\eta}, \boldsymbol{\xi})_W + \tau \int_\Sigma \operatorname{tr}_\Sigma^- \boldsymbol{\eta} \, \operatorname{tr}_\Sigma^- \boldsymbol{\xi} + \tau \int_{Q_0} \operatorname{tr}_{Q_0} \boldsymbol{\eta} \, \operatorname{tr}_{Q_0} \boldsymbol{\xi}, \tag{5.42}$$

where

$$(\boldsymbol{\eta}, \boldsymbol{\xi})_W = \int_Q \boldsymbol{\eta}^T \boldsymbol{\xi} + \tau^2 \int_Q (\mathbf{T}\boldsymbol{\eta})^T \mathbf{T}\boldsymbol{\xi}. \tag{5.43}$$

The space $C^1(\bar{Q})^{d+1}$ is dense in $V$ if and only if the only $\boldsymbol{\xi} \in V$ that satisfies

$$(\boldsymbol{\eta}, \boldsymbol{\xi})_V = 0 \quad \text{for all } \boldsymbol{\eta} \in C^1(\bar{Q})^{d+1} \tag{5.44}$$

is $\boldsymbol{\xi} = \mathbf{0}$. Let us therefore assume that $\boldsymbol{\xi} \in V$ satisfies equation (5.44) and demonstrate that $\boldsymbol{\xi} = \mathbf{0}$. In particular, equation (5.44) implies that for all $\boldsymbol{\phi} \in C_0^1(Q)^{d+1} \subset C^1(\bar{Q})^{d+1}$,

$$0 = (\boldsymbol{\phi}, \boldsymbol{\xi})_V = (\boldsymbol{\phi}, \boldsymbol{\xi})_W = \int_Q \boldsymbol{\phi}^T \boldsymbol{\xi} + \tau^2 \int_Q (\mathbf{T}\boldsymbol{\phi})^T \mathbf{T}\boldsymbol{\xi}, \tag{5.45}$$

which by the definition of weak derivatives implies that $\boldsymbol{\xi}$ satisfies the equation

$$\boldsymbol{\xi} - \tau^2 \mathbf{T}^2 \boldsymbol{\xi} = 0 \quad \text{in } Q. \tag{5.46}$$

The next step is to determine, in addition to equation (5.46), the space–time boundary conditions that $\boldsymbol{\xi}$ must satisfy. Equation (5.46) shows that $\tau \mathbf{T}\boldsymbol{\xi} \in W$ (since $\boldsymbol{\xi} \in L^2(Q)^{d+1}$). Multiplying equation (5.46) with $\boldsymbol{\eta}^T \in H^1(Q)^{d+1}$, integrating, and applying integration-by-parts formula (5.23) (in which we substitute $\boldsymbol{\xi}$ with $\tau^2 \mathbf{T}\boldsymbol{\xi}$), we obtain

$$\begin{aligned}
0 &= \int_Q \boldsymbol{\eta}^T \left(\boldsymbol{\xi} - \tau^2 \mathbf{T}^2 \boldsymbol{\xi}\right) = \int_Q \boldsymbol{\eta}^T \boldsymbol{\xi} + \tau^2 \int_Q (\mathbf{T}\boldsymbol{\eta})^T \mathbf{T}\boldsymbol{\xi} - \tau^2 \langle \mathbf{T}_\nu \mathbf{T}\boldsymbol{\xi}, \gamma\boldsymbol{\eta} \rangle_{H^{1/2}(\partial Q)^{d+1}} \\
&= (\boldsymbol{\eta}, \boldsymbol{\xi})_W - \tau^2 \langle \mathbf{T}_\nu \mathbf{T}\boldsymbol{\xi}, \gamma\boldsymbol{\eta} \rangle_{H^{1/2}(\partial Q)^{d+1}} \quad \text{for all } \boldsymbol{\eta} \in H^1(Q)^{d+1}.
\end{aligned} \tag{5.47}$$

Since $C^1(\bar{Q})^{d+1}$ is a dense subspace of $H^1(Q)^{d+1}$ and the boundary trace operator $\gamma$ on $H^1(Q)^{d+1}$ is continuous, equation (5.44) can be extended so that

$$\begin{aligned}
0 &= (\boldsymbol{\eta}, \boldsymbol{\xi})_V \\
&= (\boldsymbol{\eta}, \boldsymbol{\xi})_W + \tau \int_\Sigma \operatorname{tr}_\Sigma^- \boldsymbol{\eta} \, \operatorname{tr}_\Sigma^- \boldsymbol{\xi} + \tau \int_{Q_0} (\operatorname{tr}_{Q_0} \boldsymbol{\eta})^T \operatorname{tr}_{Q_0} \boldsymbol{\xi} \quad \text{for all } \boldsymbol{\eta} \in H^1(Q)^{d+1}.
\end{aligned} \tag{5.48}$$

Note that for $\boldsymbol{\eta} \in H^1(Q)^{d+1}$ the trace maps are given by

$$\operatorname{tr}_\Sigma^\pm \boldsymbol{\eta} = \sqrt{\frac{c_0}{2}} \left(\gamma\eta_2 \pm \boldsymbol{n} \cdot \gamma\boldsymbol{\eta}_1\right)\big|_\Sigma, \quad \operatorname{tr}_{Q_0} \boldsymbol{\eta} = \gamma\boldsymbol{\eta}\big|_{Q_0}, \quad \operatorname{tr}_{Q_\tau} \boldsymbol{\eta} = \gamma\boldsymbol{\eta}\big|_{Q_\tau}. \tag{5.49}$$

Subtracting equation (5.47) from equation (5.48), we find that

$$\tau^2 \langle \mathbf{T}_\nu \mathbf{T}\boldsymbol{\xi}, \gamma\boldsymbol{\eta}\rangle_{H^{1/2}(\partial Q)^{d+1}} = -\tau \int_\Sigma \mathrm{tr}_\Sigma^- \boldsymbol{\eta} \, \mathrm{tr}_\Sigma^- \boldsymbol{\xi} - \tau \int_{Q_0} (\mathrm{tr}_{Q_0} \boldsymbol{\eta})^T \, \mathrm{tr}_{Q_0} \boldsymbol{\xi} \quad (5.50)$$

for all $\boldsymbol{\eta} \in H^1(Q)^{d+1}$. To reveal the space–time boundary conditions satisfied by $\boldsymbol{\xi}$ on $\partial Q$ from equation (5.50), we proceed similarly as in the proof of lemma 5.7. We start with the boundary conditions on $\Sigma$. Let $\hat{\boldsymbol{\eta}}_\pm = \gamma_* \mathbf{h}_\pm \psi_* \in H^1(Q)^{d+1}$, where $\psi_* \in H^{1/2}(\partial Q)$ is the extension by zero of some $\psi \in H_{00}^{1/2}(\Sigma)$, and $\mathbf{h}_\pm$ is defined by expression (5.30). Then, by expressions (5.49), $\mathrm{tr}_\Sigma^- \hat{\boldsymbol{\eta}}_- = -\psi$, $\mathrm{tr}_\Sigma^- \hat{\boldsymbol{\eta}}_+ = 0$, and $\mathrm{tr}_{Q_0} \hat{\boldsymbol{\eta}}_\pm = 0$. Moreover, by identity (5.35) (with $\tau^2 \mathbf{T}\boldsymbol{\xi} \in W$ in place of $\boldsymbol{\xi}$), $\tau^2 \langle \mathbf{T}_\nu \mathbf{T}\boldsymbol{\xi}, \gamma\hat{\boldsymbol{\eta}}_\pm\rangle_{H^{1/2}(\partial Q)^{d+1}} = \tau^2 \langle \mathrm{tr}_\Sigma^\pm \mathbf{T}\boldsymbol{\xi}, \psi\rangle_{H_{00}^{1/2}(\Sigma)}$. Substituting these equalities into equation (5.50) with $\boldsymbol{\eta} = \hat{\boldsymbol{\eta}}_-$ and $\hat{\boldsymbol{\eta}}_+$, respectively, we find that

$$\langle \mathrm{tr}_\Sigma^- \boldsymbol{\xi} - \tau \, \mathrm{tr}_\Sigma^- \mathbf{T}\boldsymbol{\xi}, \psi\rangle_{H_{00}^{1/2}(\Sigma)} = 0 \quad \text{for all } \psi \in H_{00}^{1/2}(\Sigma), \quad (5.51)$$

$$\langle \mathrm{tr}_\Sigma^+ \mathbf{T}\boldsymbol{\xi}, \psi\rangle_{H_{00}^{1/2}(\Sigma)} = 0 \quad \text{for all } \psi \in H_{00}^{1/2}(\Sigma). \quad (5.52)$$

That is,

$$\mathrm{tr}_\Sigma^-(\boldsymbol{\xi} - \tau \mathbf{T}\boldsymbol{\xi}) = 0 \quad \text{on } \Sigma, \quad (5.53)$$

$$\mathrm{tr}_\Sigma^+ \tau \mathbf{T}\boldsymbol{\xi} = 0 \quad \text{on } \Sigma. \quad (5.54)$$

Analogously, by letting, in equation (5.50), $\boldsymbol{\eta} = \gamma_* \boldsymbol{\psi}_* \in H^1(Q)^{d+1}$, where $\boldsymbol{\psi}_* \in H^{1/2}(\partial Q)^{d+1}$ is the extension by zeros of some $\boldsymbol{\psi} \in H_{00}^{1/2}(Q_0)^{d+1}$, and recalling identity (5.33) (with $\mathbf{T}\boldsymbol{\xi}$ in place of $\boldsymbol{\xi}$), we find that

$$\mathrm{tr}_{Q_0}(\boldsymbol{\xi} - \tau \mathbf{T}\boldsymbol{\xi}) = \mathbf{0} \quad \text{on } Q_0, \quad (5.55)$$

while $\hat{\boldsymbol{\eta}} = \gamma_* \boldsymbol{\psi}_* \in H^1(Q)^{d+1}$, where $\boldsymbol{\psi}_* \in H^{1/2}(\partial Q)^{d+1}$ is the extension by zero of some $\boldsymbol{\psi} \in H_{00}^{1/2}(Q_\tau)^{d+1}$, in equation (5.50) implies (recall identity (5.34)) that

$$\mathrm{tr}_{Q_\tau} \tau \mathbf{T}\boldsymbol{\xi} = \mathbf{0} \quad \text{on } Q_\tau. \quad (5.56)$$

Let $\mathrm{id}_W$ denote the identity operator on $W$, that is, $\mathrm{id}_W \boldsymbol{\psi} = \boldsymbol{\psi}$ for any $\boldsymbol{\psi} \in W$. Note that

$$(\mathrm{id}_W + \tau \mathbf{T})(\mathrm{id}_W - \tau \mathbf{T})\boldsymbol{\xi} = \boldsymbol{\xi} - \tau^2 \mathbf{T}\boldsymbol{\xi}, \quad (5.57)$$

so by introducing $\boldsymbol{\psi} = (\mathrm{id}_W - \tau \mathbf{T})\boldsymbol{\xi} = \boldsymbol{\xi} - \tau \mathbf{T}\boldsymbol{\xi} \in W$, we may reformulate the 2nd-order problem formed by equation (5.46) and conditions (5.53), (5.54), (5.55), and (5.56) as the coupled 1st-order system

$$\boldsymbol{\psi} + \tau \mathbf{T}\boldsymbol{\psi} = \mathbf{0} \quad \text{in } Q, \quad (5.58a)$$

$$\mathrm{tr}_\Sigma^- \boldsymbol{\psi} = \mathbf{0} \quad \text{on } \Sigma, \quad (5.58b)$$

$$\mathrm{tr}_{Q_0} \boldsymbol{\psi} = \mathbf{0} \quad \text{on } Q_0, \quad (5.58c)$$

and

$$\boldsymbol{\xi} - \tau \mathbf{T}\boldsymbol{\xi} = \boldsymbol{\psi} \quad \text{in } Q, \quad (5.59a)$$

$$\mathrm{tr}_\Sigma^+ \tau \mathbf{T}\boldsymbol{\xi} = \mathbf{0} \quad \text{on } \Sigma, \quad (5.59b)$$

$$\mathrm{tr}_{Q_\tau} \tau \mathbf{T}\boldsymbol{\xi} = \mathbf{0} \quad \text{on } Q_T. \quad (5.59c)$$

Recall that $\boldsymbol{\xi}$, $\tau \mathbf{T}\boldsymbol{\xi}$, and $\boldsymbol{\psi}$ belong to $W$. Thus, by applying $\mathrm{tr}_\Sigma^+$ to equation (5.59a) and invoking boundary condition (5.59b), we find that

$$\mathrm{tr}_\Sigma^+ \boldsymbol{\xi} = \mathrm{tr}_\Sigma^+ (\tau \mathbf{T}\boldsymbol{\xi} + \boldsymbol{\psi}) = \mathrm{tr}_\Sigma^+ \boldsymbol{\psi} \quad \text{on } \Sigma. \quad (5.60)$$

Analogously, by applying $\mathrm{tr}_{Q_\tau}$ to equation (5.59a) and invoking final condition (5.59c), we find that

$$\mathrm{tr}_{Q_\tau}\,\boldsymbol{\xi} = \mathrm{tr}_{Q_\tau}(\tau \mathbf{T}\boldsymbol{\xi} + \boldsymbol{\psi}) = \mathrm{tr}_{Q_\tau}\,\boldsymbol{\psi} \quad \text{on } Q_\tau. \tag{5.61}$$

Thus, system (5.59) translates to the system

$$\boldsymbol{\xi} - \tau \mathbf{T}\boldsymbol{\xi} = \boldsymbol{\psi} \quad \text{in } Q, \tag{5.62a}$$
$$\mathrm{tr}_\Sigma^+\,\boldsymbol{\xi} = \mathrm{tr}_\Sigma^+\,\boldsymbol{\psi} \quad \text{on } \Sigma, \tag{5.62b}$$
$$\mathrm{tr}_{Q_\tau}\,\boldsymbol{\xi} = \mathrm{tr}_{Q_\tau}\,\boldsymbol{\psi} \quad \text{on } Q_T. \tag{5.62c}$$

Note that $\boldsymbol{\psi} \in V_0$, so by theorem 5.9 there exists $(\boldsymbol{\psi}_k)_{k \in \mathbb{Z}^+}$, where $\boldsymbol{\psi}_k = [\psi_{k,1}, \psi_{k,2}]^T \in C^1(\overline{Q})^{d+1}$ such that $\mathrm{tr}_\Sigma^-\,\boldsymbol{\psi}_k = 0$, $\mathrm{tr}_{Q_0}\,\boldsymbol{\psi}_k = 0$, and $\|\boldsymbol{\psi}_k - \boldsymbol{\psi}\|_W \to 0$ when $k \to \infty$. Integration-by-parts formula (5.21) yields

$$\begin{aligned}
\int_Q \boldsymbol{\psi}_k^T \tau \mathbf{T}\boldsymbol{\psi}_k &= \frac{\tau}{2}\int_{\partial Q} \boldsymbol{\psi}_k^T \mathbf{T}_\nu \boldsymbol{\psi}_k \\
&= -\frac{\tau}{2}\int_{Q_0} |\mathrm{tr}_{Q_0}\,\boldsymbol{\psi}_k|^2 + \frac{\tau}{2}\int_{Q_\tau} |\mathrm{tr}_{Q_\tau}\,\boldsymbol{\psi}_k|^2 + \tau \int_\Sigma c_0\, \boldsymbol{n}\cdot \boldsymbol{\psi}_{k,1}\psi_{k,2} \\
&= -\frac{\tau}{2}\int_{Q_0} \underbrace{|\mathrm{tr}_{Q_0}\,\boldsymbol{\psi}_k|^2}_{=0} + \frac{\tau}{2}\int_{Q_\tau} |\mathrm{tr}_{Q_\tau}\,\boldsymbol{\psi}_k|^2 + \frac{\tau}{2}\int_\Sigma \left(|\mathrm{tr}_\Sigma^+\,\boldsymbol{\psi}_k|^2 - \underbrace{|\mathrm{tr}_\Sigma^-\,\boldsymbol{\psi}_k|^2}_{=0}\right) \\
&\geq 0.
\end{aligned} \tag{5.63}$$

where identity

$$2ab = \frac{(a+b)^2 - (a-b)^2}{2} \tag{5.64}$$

with $a = \sqrt{c_0/2}\,\boldsymbol{n}\cdot\boldsymbol{\psi}_{k,1}$ and $b = \sqrt{c_0/2}\,\psi_{k,2}$ is used in the third equality. Thus,

$$\|\boldsymbol{\psi}_k\|_{L^2(Q)^{d+1}}^2 = \int_Q |\boldsymbol{\psi}_k|^2 \leq \int_Q \boldsymbol{\psi}_k^T(\boldsymbol{\psi}_k + \tau \mathbf{T}\boldsymbol{\psi}_k) \leq C\|\boldsymbol{\psi}_k + \tau\mathbf{T}\boldsymbol{\psi}_k\|_{L^2(Q)^{d+1}}, \tag{5.65}$$

where the bound (5.63) is used in the first inequality, and second inequality follows from the Cauchy–Schwarz inequality and by choosing a constant $C$ such that $\|\boldsymbol{\psi}_k\|_{L^2(Q)^{d+1}} \leq C$ for all $k \in \mathbb{Z}^+$. Passing to the limit in estimate (5.65), recalling that $\boldsymbol{\psi}$ satisfies equation (5.58a), demonstrates that $\boldsymbol{\psi} = \mathbf{0}$. Therefore, problem (5.62) reads

$$\boldsymbol{\xi} - \tau \mathbf{T}\boldsymbol{\xi} = 0 \quad \text{in } Q, \tag{5.66a}$$
$$\mathrm{tr}_\Sigma^+\,\boldsymbol{\xi} = 0 \quad \text{on } \Sigma, \tag{5.66b}$$
$$\mathrm{tr}_{Q_\tau}\,\boldsymbol{\xi} = 0 \quad \text{on } Q_T. \tag{5.66c}$$

Proceeding similarly as for $\boldsymbol{\psi}$ above, we note that $\boldsymbol{\xi} \in \widetilde{V}_0$, so by theorem 5.9 there exists $(\boldsymbol{\xi}_k)_{k\in\mathbb{Z}^+}$, where $\boldsymbol{\xi}_k = [\xi_{k,1}, \xi_{k,2}]^T \in C^1(\overline{Q})^{d+1}$ such that $\mathrm{tr}_\Sigma^+\,\boldsymbol{\xi}_k = 0$, $\mathrm{tr}_{Q_\tau}\,\boldsymbol{\xi}_k = 0$, and $\|\boldsymbol{\xi}_k - \boldsymbol{\xi}\|_W \to 0$ when $k \to \infty$. Integration-by-parts formula (5.21) yields

$$\begin{aligned}
-\int_Q \boldsymbol{\xi}_k^T \tau \mathbf{T}\boldsymbol{\xi}_k &= -\frac{\tau}{2}\int_{\partial Q} \boldsymbol{\xi}_k^T \mathbf{T}_\nu \boldsymbol{\xi}_k \\
&= \frac{\tau}{2}\int_{Q_0} |\mathrm{tr}_{Q_0}\,\boldsymbol{\xi}_k|^2 - \frac{\tau}{2}\int_{Q_\tau} \underbrace{|\mathrm{tr}_{Q_\tau}\,\boldsymbol{\xi}_k|^2}_{=0} - \frac{\tau}{2}\int_\Sigma \left(\underbrace{|\mathrm{tr}_\Sigma^+\,\boldsymbol{\xi}_k|^2}_{=0} - |\mathrm{tr}_\Sigma^-\,\boldsymbol{\xi}_k|^2\right) \\
&\geq 0,
\end{aligned} \tag{5.67}$$

where also here identity (5.64) is used to obtain the last term after the second equality. Thus, similarly as in expression (5.65), also here we arrive at the analogous bound

$$\|\boldsymbol{\xi}_k\|_{L^2(Q)^{d+1}}^2 = \int_Q |\boldsymbol{\xi}_k|^2 \leq \int_Q \boldsymbol{\xi}_k^T(\boldsymbol{\xi}_k - \tau \mathbf{T}\boldsymbol{\xi}_k) \leq C\|\boldsymbol{\xi}_k - \tau\mathbf{T}\boldsymbol{\xi}_k\|_{L^2(Q)^{d+1}}, \tag{5.68}$$

where constant $C$ is such that $\|\boldsymbol{\xi}_k\|_{L^2(Q)^{d+1}} \leq C$ for all $k \in \mathbb{Z}^+$. Passing to the limit in estimate (5.68), recalling that $\boldsymbol{\xi}$ satisfies equation (5.66a), finally demonstrates that $\boldsymbol{\xi} = \mathbf{0}$. $\square$

Although only trace operators $\operatorname{tr}_{Q_0}$ and $\operatorname{tr}_{\Sigma}^-$ are involved in the definition of the space $V$, it turns out that the remaining trace operators also map continuously into $L^2$ spaces.

**Lemma 5.11.** *The trace operators $\operatorname{tr}_{Q_\tau}$ and $\operatorname{tr}_{\Sigma}^+$, defined in expression (5.16) for functions in $C^1(\bar{Q})^{d+1}$, extend continuously to $\mathscr{L}(V, L^2(Q_\tau)^{d+1})$ and $\mathscr{L}(V, L^2(\Sigma))$, respectively.*

*Proof.* Integration-by-parts formula (5.21) implies that for $\boldsymbol{\xi} \in C^1(\bar{Q})^{d+1}$,

$$
\begin{aligned}
\int_Q \boldsymbol{\xi}^T \mathbf{T}\boldsymbol{\xi} &= \frac{1}{2} \int_{Q_\tau} |\operatorname{tr}_{Q_\tau} \boldsymbol{\xi}|^2 - \frac{1}{2} \int_{Q_0} |\operatorname{tr}_{Q_0} \boldsymbol{\xi}|^2 + c_0 \int_\Sigma \boldsymbol{n} \cdot \boldsymbol{\xi}_1 \xi_2 \\
&= \frac{1}{2} \int_{Q_\tau} |\operatorname{tr}_{Q_\tau} \boldsymbol{\xi}|^2 - \frac{1}{2} \int_{Q_0} |\operatorname{tr}_{Q_0} \boldsymbol{\xi}|^2 + \frac{c_0}{4} \int_\Sigma \left[ (\xi_2 + \boldsymbol{n} \cdot \boldsymbol{\xi}_1)^2 - (\xi_2 - \boldsymbol{n} \cdot \boldsymbol{\xi}_1)^2 \right] \\
&= \frac{1}{2} \int_{Q_\tau} |\operatorname{tr}_{Q_\tau} \boldsymbol{\xi}|^2 - \frac{1}{2} \int_{Q_0} |\operatorname{tr}_{Q_0} \boldsymbol{\xi}|^2 + \frac{1}{2} \int_\Sigma \left[ (\operatorname{tr}_{\Sigma}^+ \boldsymbol{\xi})^2 - (\operatorname{tr}_{\Sigma}^- \boldsymbol{\xi})^2 \right],
\end{aligned}
\quad (5.69)
$$

again using identity (5.64) in the second equality, which implies that

$$
\begin{aligned}
\tau \int_{Q_\tau} |\operatorname{tr}_{Q_\tau} \boldsymbol{\xi}|^2 + \tau \int_\Sigma (\operatorname{tr}_{\Sigma}^+ \boldsymbol{\xi})^2 &= \tau \int_{Q_0} |\operatorname{tr}_{Q_0} \boldsymbol{\xi}|^2 + \tau \int_\Sigma (\operatorname{tr}_{\Sigma}^- \boldsymbol{\xi})^2 + 2\tau \int_Q \boldsymbol{\xi}^T \mathbf{T}\boldsymbol{\xi} \\
&\leq \tau \int_{Q_0} |\operatorname{tr}_{Q_0} \boldsymbol{\xi}|^2 + \tau \int_\Sigma (\operatorname{tr}_{\Sigma}^- \boldsymbol{\xi})^2 + \int_Q |\boldsymbol{\xi}|^2 + \tau^2 \int_Q |\mathbf{T}\boldsymbol{\xi}|^2 = \|\boldsymbol{\xi}\|_V^2,
\end{aligned}
\quad (5.70)
$$

by definition (5.40) of the norm on $V$, from which the conclusion follows by density theorem 5.10. $\square$

By lemma 5.11, we conclude that bilinear forms $a$ and $a^*$ are well defined for $\boldsymbol{\xi} \in V$, and we are ready to show well-posedness for the variational problem in standard form:

$$
\begin{aligned}
&\text{Find } \boldsymbol{\xi} \in V \text{ such that} \\
&a(\hat{\boldsymbol{\eta}}, \boldsymbol{\xi}) = l(\hat{\boldsymbol{\eta}}) \qquad \forall \hat{\boldsymbol{\eta}} \in L.
\end{aligned}
\quad (5.71)
$$

The first step to acquire the inf–sup condition is the following bounds.

**Lemma 5.12.** *For each $\boldsymbol{\xi} \in V$,*

$$
\sup_{\substack{\hat{\boldsymbol{\eta}} \in L \\ \hat{\boldsymbol{\eta}} \neq 0}} \frac{a(\hat{\boldsymbol{\eta}}, \boldsymbol{\xi})}{\|\hat{\boldsymbol{\eta}}\|_L} \geq \beta \|\hat{\boldsymbol{\xi}}\|_L, \quad (5.72a)
$$

$$
\sup_{\substack{\hat{\boldsymbol{\eta}} \in L^* \\ \hat{\boldsymbol{\eta}} \neq 0}} \frac{a^*(\hat{\boldsymbol{\eta}}, \boldsymbol{\xi})}{\|\hat{\boldsymbol{\eta}}\|_{L^*}} \geq \beta \|\hat{\boldsymbol{\xi}}^*\|_{L^*}, \quad (5.72b)
$$

*where*

$$
\beta = \frac{1}{4\tau \mathrm{e}}(1 - \alpha_M),
$$
$$
\hat{\boldsymbol{\xi}} = (\boldsymbol{\xi}, \operatorname{tr}_{\Sigma}^- \boldsymbol{\xi}, \operatorname{tr}_{Q_0} \boldsymbol{\xi}), \qquad \hat{\boldsymbol{\xi}}^* = (\boldsymbol{\xi}, \operatorname{tr}_{\Sigma}^+ \boldsymbol{\xi}, \operatorname{tr}_{Q_\tau} \boldsymbol{\xi})
$$
(5.73)

*Proof.* Due to density theorem 5.10, it is enough to show the inequalities for $\boldsymbol{\xi} \in C^1(\bar{Q})^{d+1}$. Since the statements are immediate for $\boldsymbol{\xi} = \mathbf{0}$, let $\boldsymbol{\xi} \in C^1(\bar{Q})^{d+1}$ be nonzero, define

$$
\tilde{\boldsymbol{\eta}} = \left( \mathrm{e}^{-t/\tau} \boldsymbol{\xi}, \mathrm{e}^{-t/\tau} \operatorname{tr}_{\Sigma}^- \boldsymbol{\xi}, \operatorname{tr}_{Q_0} \boldsymbol{\xi} \right), \quad (5.74)
$$

and note that

$$
\|\tilde{\boldsymbol{\eta}}\|_L \geq \mathrm{e}^{-1} \|\hat{\boldsymbol{\xi}}\|_L. \quad (5.75)
$$

Choosing $\hat{\boldsymbol{\eta}} = \widetilde{\boldsymbol{\eta}}$ in definition (5.17a), we find that

$$a(\widetilde{\boldsymbol{\eta}}, \boldsymbol{\xi}) = \int_Q e^{-t/\tau} \boldsymbol{\xi}^T \mathbf{T}\boldsymbol{\xi} + \int_\Sigma e^{-t/\tau} \operatorname{tr}_\Sigma^- \boldsymbol{\xi} \left( \operatorname{tr}_\Sigma^- \boldsymbol{\xi} - \alpha \operatorname{tr}_\Sigma^+ \boldsymbol{\xi} \right) + \int_{Q_0} |\operatorname{tr}_{Q_0} \boldsymbol{\xi}|^2. \tag{5.76}$$

Applying integration-by-parts formula (5.21), the first term in expression (5.76) can be written

$$\int_Q e^{-t/\tau} \boldsymbol{\xi}^T \mathbf{T}\boldsymbol{\xi} = \int_{\partial Q} e^{-t/\tau} \boldsymbol{\xi}^T \mathbf{T}_\nu \boldsymbol{\xi} - \int_Q \boldsymbol{\xi}^T \mathbf{T}(e^{-t/\tau} \boldsymbol{\xi}), \tag{5.77}$$

from which it follows, after substituting the identity

$$\mathbf{T}(e^{-t/\tau} \boldsymbol{\xi}) = e^{-t/\tau} \mathbf{T}\boldsymbol{\xi} - \frac{e^{-t/\tau}}{\tau} \boldsymbol{\xi}, \tag{5.78}$$

that

$$\begin{aligned}
\int_Q e^{-t/\tau} \boldsymbol{\xi}^T \mathbf{T}\boldsymbol{\xi} &= \frac{1}{2} \int_{\partial Q} e^{-t/\tau} \boldsymbol{\xi}^T \mathbf{T}_\nu \boldsymbol{\xi} + \frac{1}{2\tau} \int_Q e^{-t/\tau} |\boldsymbol{\xi}|^2 \\
&= \frac{1}{2e} \int_{Q_\tau} |\operatorname{tr}_{Q_\tau} \boldsymbol{\xi}|^2 - \frac{1}{2} \int_{Q_0} |\operatorname{tr}_{Q_0} \boldsymbol{\xi}|^2 + \frac{1}{2} \int_\Sigma e^{-t/\tau} c_0 \, \boldsymbol{n} \cdot \boldsymbol{\xi}_1 \xi_2 + \frac{1}{2\tau} \int_Q e^{-t/\tau} |\boldsymbol{\xi}|^2,
\end{aligned} \tag{5.79}$$

where definition (5.22) is used in the second equality. Substituting equality (5.79) into expression (5.76), we find that

$$\begin{aligned}
a(\widetilde{\boldsymbol{\eta}}, \boldsymbol{\xi}) &= \frac{1}{2\tau} \int_Q e^{-t/\tau} |\boldsymbol{\xi}|^2 + \frac{1}{2e} \int_{Q_\tau} |\operatorname{tr}_{Q_\tau} \boldsymbol{\xi}|^2 + \frac{1}{2} \int_{Q_0} |\operatorname{tr}_{Q_0} \boldsymbol{\xi}|^2 \\
&\quad + \frac{1}{2} \int_\Sigma e^{-t/\tau} c_0 \, \boldsymbol{n} \cdot \boldsymbol{\xi}_1 \xi_2 + \frac{1}{2} \int_\Sigma e^{-t/\tau} \operatorname{tr}_\Sigma^- \boldsymbol{\xi} \left( \operatorname{tr}_\Sigma^- \boldsymbol{\xi} - \alpha \operatorname{tr}_\Sigma^+ \boldsymbol{\xi} \right).
\end{aligned} \tag{5.80}$$

Now, since

$$c_0 \, \boldsymbol{n} \cdot \boldsymbol{\xi}_1 \xi_2 = \frac{c_0}{4} \left[ (\xi_2 + \boldsymbol{n} \cdot \boldsymbol{\xi}_1)^2 - (\xi_2 - \boldsymbol{n} \cdot \boldsymbol{\xi}_1)^2 \right] = \frac{1}{2} \left[ (\operatorname{tr}_\Sigma^+ \boldsymbol{\xi})^2 - (\operatorname{tr}_\Sigma^- \boldsymbol{\xi})^2 \right], \tag{5.81}$$

the last two integrals in expression (5.80) can be written as

$$\begin{aligned}
\frac{1}{2} \int_\Sigma e^{-t/\tau} c_0 \, \boldsymbol{n} \cdot \boldsymbol{\xi}_1 \xi_2 &+ \frac{1}{2} \int_\Sigma e^{-t/\tau} \operatorname{tr}_\Sigma^- \boldsymbol{\xi} \left( \operatorname{tr}_\Sigma^- \boldsymbol{\xi} - \alpha \operatorname{tr}_\Sigma^+ \boldsymbol{\xi} \right) \\
&= \frac{1}{4} \int_\Sigma e^{-t/\tau} \left[ (\operatorname{tr}_\Sigma^+ \boldsymbol{\xi})^2 - (\operatorname{tr}_\Sigma^- \boldsymbol{\xi})^2 \right] + \frac{1}{2} \int_\Sigma e^{-t/\tau} \operatorname{tr}_\Sigma^- \boldsymbol{\xi} \left( \operatorname{tr}_\Sigma^- \boldsymbol{\xi} - \alpha \operatorname{tr}_\Sigma^+ \boldsymbol{\xi} \right). \\
&\geq \frac{1 - \alpha_M}{4} \int_\Sigma e^{-t/\tau} \left( \operatorname{tr}_\Sigma^- \boldsymbol{\xi} \right)^2 \geq \frac{1 - \alpha_M}{4} \int_\Sigma \left( e^{-t/\tau} \operatorname{tr}_\Sigma^- \boldsymbol{\xi} \right)^2,
\end{aligned} \tag{5.82}$$

where the first inequality follows from setting $a = \operatorname{tr}_\Sigma^- \boldsymbol{\xi}$, $b = \operatorname{tr}_\Sigma^+ \boldsymbol{\xi}$ and observing that

$$\begin{aligned}
\frac{1}{4}(b^2 - a^2) + \frac{1}{2} a(a - \alpha b) &= \frac{1}{4}(a^2 + b^2 - 2\alpha ab) \geq \frac{1}{4}(a^2 + b^2 - 2\alpha_M |ab|) \\
&\geq \frac{1}{4}(a^2 + b^2 - \alpha_M(a^2 + b^2)) \geq \frac{1 - \alpha_M}{4} a^2.
\end{aligned} \tag{5.83}$$

Substituting inequality (5.82) into equality (5.80), we arrive at the bound

$$\begin{aligned}
a(\widetilde{\boldsymbol{\eta}}, \boldsymbol{\xi}) &\geq \frac{1}{2\tau} \int_Q |e^{-t/\tau} \boldsymbol{\xi}|^2 + \frac{1 - \alpha_M}{4} \int_\Sigma \left( e^{-t/\tau} \operatorname{tr}_\Sigma^- \boldsymbol{\xi} \right)^2 + \frac{1}{2} \int_{Q_0} |\operatorname{tr}_{Q_0} \boldsymbol{\xi}|^2 \\
&\geq \frac{1}{4\tau}(1 - \alpha_M) \|\widetilde{\boldsymbol{\eta}}\|_L^2 \geq \beta \|\widetilde{\boldsymbol{\eta}}\|_L \|\hat{\boldsymbol{\xi}}\|_L,
\end{aligned} \tag{5.84}$$

where the last inequality follows from bound (5.75). Dividing by $\|\widetilde{\boldsymbol{\eta}}\|_L$ and taking supremum yields inequality (5.72a).

Inequality (5.72b) is shown analogously. □

*Remark* 5.13. The use of the exponential weighting in time for the test functions, introduced in expression (5.74), is crucial to obtain the "$L$-coercivity" property (5.72a), and compensates for the lack of property (1.7) in this example. An alternative would be to employ an exponentially weighted Hilbert space in time, as done by Franz et al. [10].

With the help of lemma 5.12, the inf–sup condition is straightforward to achieve.

**Lemma 5.14.** *For each $\boldsymbol{\xi} \in V$,*
$$\sup_{\substack{\hat{\boldsymbol{\eta}} \in L \\ \hat{\boldsymbol{\eta}} \neq 0}} \frac{a(\hat{\boldsymbol{\eta}}, \boldsymbol{\xi})}{\|\hat{\boldsymbol{\eta}}\|_L} \geq \frac{\beta}{2} \|\boldsymbol{\xi}\|_V. \tag{5.85}$$

*Proof.* For $\mathbf{T}\boldsymbol{\xi} = \mathbf{0}$, the conclusion follows from lemma 5.12. Thus, assume $\mathbf{T}\boldsymbol{\xi} \neq \mathbf{0}$ and let $\widetilde{\boldsymbol{\eta}} = (\mathbf{T}\boldsymbol{\xi}, 0, 0)$. Then
$$\sup_{\substack{\hat{\boldsymbol{\eta}} \in L \\ \hat{\boldsymbol{\eta}} \neq 0}} \frac{a(\hat{\boldsymbol{\eta}}, \boldsymbol{\xi})}{\|\hat{\boldsymbol{\eta}}\|_L} \geq \frac{a(\widetilde{\boldsymbol{\eta}}, \boldsymbol{\xi})}{\|\widetilde{\boldsymbol{\eta}}\|_L} = \|\mathbf{T}\boldsymbol{\xi}\|_{L^2(Q)^{d+1}}. \tag{5.86}$$

Writing the norm on $V$ as in expression (5.41), we reach the conclusion from the estimates
$$\|\boldsymbol{\xi}\|_V = \left(\|(\boldsymbol{\xi}, \operatorname{tr}_\Sigma^- \boldsymbol{\xi}, \operatorname{tr}_{Q_0} \boldsymbol{\xi})\|_L^2 + \tau^2 \|\mathbf{T}\boldsymbol{\xi}\|_{L^2(Q)^{d+1}}^2\right)^{1/2}$$
$$\leq \|(\boldsymbol{\xi}, \operatorname{tr}_\Sigma^- \boldsymbol{\xi}, \operatorname{tr}_{Q_0} \boldsymbol{\xi})\|_L + \tau \|\mathbf{T}\boldsymbol{\xi}\|_{L^2(Q)^{d+1}} \leq \left(\frac{1}{\beta} + \tau\right) \sup_{\substack{\hat{\boldsymbol{\eta}} \in L \\ \hat{\boldsymbol{\eta}} \neq 0}} \frac{a(\hat{\boldsymbol{\eta}}, \boldsymbol{\xi})}{\|\hat{\boldsymbol{\eta}}\|_L} \leq \frac{2}{\beta} \frac{a(\hat{\boldsymbol{\eta}}, \boldsymbol{\xi})}{\|\hat{\boldsymbol{\eta}}\|_L}, \tag{5.87}$$

where the second inequality follows from lemma 5.12 and bound (5.86), and the third since $\tau\beta < 1$. □

Surjectivity is shown in a manner very similar to the previous two examples.

**Lemma 5.15.** *If $\hat{\boldsymbol{\eta}} \in L$ such that*
$$a(\hat{\boldsymbol{\eta}}, \boldsymbol{\xi}) = 0 \qquad \forall \boldsymbol{\xi} \in V, \tag{5.88}$$

*then $\hat{\boldsymbol{\eta}} = \mathbf{0}$.*

*Proof.* Let $\hat{\boldsymbol{\eta}} = (\boldsymbol{\eta}, \eta_\Sigma, \boldsymbol{\eta}_s) \in L$ satisfy
$$a(\hat{\boldsymbol{\eta}}, \boldsymbol{\xi}) = \int_Q \boldsymbol{\eta}^T \mathbf{T}\boldsymbol{\xi} + \int_\Sigma \eta_\Sigma (\operatorname{tr}_\Sigma^- \boldsymbol{\xi} - \alpha \operatorname{tr}_\Sigma^+ \boldsymbol{\xi}) + \int_{Q_0} \boldsymbol{\eta}_s^T (\operatorname{tr}_{Q_0} \boldsymbol{\xi}) = 0 \qquad \forall \boldsymbol{\xi} \in V. \tag{5.89}$$

We will show that all components of the tuple $\hat{\boldsymbol{\eta}}$ then vanish. The strategy is to choose various subspaces of $V$ for $\boldsymbol{\xi}$ in equation (5.89) in order to uncover information of $\hat{\boldsymbol{\eta}}$.

First, from equation (5.89) it follows that
$$a(\hat{\boldsymbol{\eta}}, \boldsymbol{\phi}) = \int_Q \boldsymbol{\eta}^T \mathbf{T}\boldsymbol{\phi} = 0 \qquad \forall \boldsymbol{\phi} \in C_0^1(Q)^{d+1}, \tag{5.90}$$

and thus that
$$\widetilde{\mathbf{T}}\boldsymbol{\eta} = -\mathbf{T}\boldsymbol{\eta} = \mathbf{0} \tag{5.91}$$

by the definition of weak derivative. Hence, trivially, $\boldsymbol{\eta} \in W$. Applying integration-by-parts-formula (5.23), with $\boldsymbol{\xi} \in H^1(Q)^{d+1} \subset V$ and $\boldsymbol{\eta} \in W$, to equation (5.89), we obtain
$$a(\hat{\boldsymbol{\eta}}, \boldsymbol{\xi}) = -\int_Q \boldsymbol{\xi}^T \mathbf{T}\boldsymbol{\eta} + \langle \mathbf{T}_\nu \boldsymbol{\eta}, \gamma \boldsymbol{\xi} \rangle_{H^{1/2}(\partial Q)^{d+1}} + \int_\Sigma \eta_\Sigma (\operatorname{tr}_\Sigma^- \boldsymbol{\xi} - \alpha \operatorname{tr}_\Sigma^+ \boldsymbol{\xi}) + \int_{Q_0} \boldsymbol{\eta}_s^T (\operatorname{tr}_{Q_0} \boldsymbol{\xi})$$
$$= \langle \mathbf{T}_\nu \boldsymbol{\eta}, \gamma \boldsymbol{\xi} \rangle_{H^{1/2}(\partial Q)^{d+1}} + \int_\Sigma \eta_\Sigma (\operatorname{tr}_\Sigma^- \boldsymbol{\xi} - \alpha \operatorname{tr}_\Sigma^+ \boldsymbol{\xi}) + \int_{Q_0} \boldsymbol{\eta}_s^T (\operatorname{tr}_{Q_0} \boldsymbol{\xi}) = 0. \tag{5.92}$$

As before, we denote by $\gamma$ the trace operator from $H^1(Q)$ onto $H^{1/2}(\partial Q)$ and by $\gamma_*$ a continuous right inverse of $\gamma$. Let $\boldsymbol{\xi}_\tau \in H_{00}^{1/2}(Q_\tau)^{d+1}$, and let $\boldsymbol{\xi}_{\tau*} \in H^{1/2}(\partial Q)^{d+1}$ be its extension by zero. Then $\gamma_* \boldsymbol{\xi}_{\tau*} \in H^1(Q)^{d+1}$ and, since $\boldsymbol{\xi}_{\tau*}$ vanishes on $\Sigma \cup Q_0$, by equation (5.92) and identity (5.34), we find
$$a(\hat{\boldsymbol{\eta}}, \gamma_* \boldsymbol{\xi}_{\tau*}) = \langle \mathbf{T}_\nu \boldsymbol{\eta}, \boldsymbol{\xi}_{\tau*} \rangle_{H^{1/2}(\partial Q)^{d+1}} = \langle \operatorname{tr}_{Q_\tau} \boldsymbol{\eta}, \boldsymbol{\xi}_\tau \rangle_{H_{00}^{1/2}(Q_\tau)^{d+1}} = 0 \tag{5.93}$$

for any $\boldsymbol{\xi}_\tau \in H_{00}^{1/2}(Q_\tau)^{d+1}$, that is,
$$\operatorname{tr}_{Q_\tau} \boldsymbol{\eta} = 0. \tag{5.94}$$

An analogous procedure (recall identity (5.33)) shows that, for any $\boldsymbol{\xi}_0 \in H_{00}^{1/2}(Q_0)^{d+1}$,

$$\begin{aligned}a(\hat{\boldsymbol{\eta}}, -\gamma_* \boldsymbol{\xi}_{0*}) &= \langle \mathbf{T}_\nu \boldsymbol{\eta}, -\boldsymbol{\xi}_{0*}\rangle_{H^{1/2}(\partial Q)^{d+1}} - \int_{Q_0} \boldsymbol{\eta}_s^T \boldsymbol{\xi}_0 \\ &= \langle \operatorname{tr}_{Q_0} \boldsymbol{\eta} - \boldsymbol{\eta}_s, \boldsymbol{\xi}_0 \rangle_{H_{00}^{1/2}(Q_0)^{d+1}} = 0,\end{aligned} \tag{5.95}$$

and thus that
$$\operatorname{tr}_{Q_0} \boldsymbol{\eta} - \boldsymbol{\eta}_s = 0, \tag{5.96}$$

which means, in particular, that $\operatorname{tr}_{Q_0} \boldsymbol{\eta} \in L^2(Q_0)^{d+1}$. Finally, let $\psi \in H_{00}^{1/2}(\Sigma)$, $\psi_* \in H^{1/2}(\partial Q)$ its extension by zero, and define $\boldsymbol{\xi}_\Sigma^\pm = \gamma_* \psi_* \mathbf{h}_\pm$, where $\mathbf{h}_\pm \in \operatorname{Lip}(\partial Q)^{d+1}$ is defined by expression (5.30). By lemma 5.6, we have that $\boldsymbol{\xi}_\Sigma^\pm \in H^1(Q)^{d+1}$. Moreover, expressions (5.49) reveal that $\operatorname{tr}_{Q_0} \boldsymbol{\xi}_\Sigma^\pm = \operatorname{tr}_{Q_\tau} \boldsymbol{\xi}_\Sigma^\pm = \operatorname{tr}_\Sigma^\mp \boldsymbol{\xi}_\Sigma^\pm = 0$ and $\operatorname{tr}_\Sigma^\pm \boldsymbol{\xi}_\Sigma^\pm = \pm \psi$. Thus, by equation (5.92) and identity (5.35), for any $\psi \in H_{00}^{1/2}(\Sigma)$,

$$\begin{aligned}a(\hat{\boldsymbol{\eta}}, \boldsymbol{\xi}_\Sigma^-) &= \langle \mathbf{T}_\nu \boldsymbol{\eta}, \gamma \boldsymbol{\xi}_\Sigma^-\rangle_{H^{1/2}(\partial Q)^{d+1}} + \int_\Sigma \eta_\Sigma \operatorname{tr}_\Sigma^- \boldsymbol{\xi}_\Sigma^- \\ &= \langle \operatorname{tr}_\Sigma^- \boldsymbol{\eta} - \eta_\Sigma, \psi\rangle_{H_{00}^{1/2}(\Sigma)} = 0,\end{aligned} \tag{5.97}$$

and

$$\begin{aligned}a(\hat{\boldsymbol{\eta}}, \boldsymbol{\xi}_\Sigma^+) &= \langle \mathbf{T}_\nu \boldsymbol{\eta}, \gamma \boldsymbol{\xi}_\Sigma^+\rangle_{H^{1/2}(\partial Q)^{d+1}} - \int_\Sigma \eta_\Sigma \alpha \operatorname{tr}_\Sigma^+ \boldsymbol{\xi}_\Sigma^+ \\ &= \langle \operatorname{tr}_\Sigma^+ \boldsymbol{\eta} - \alpha \eta_\Sigma, \psi\rangle_{H_{00}^{1/2}(\Sigma)} = 0,\end{aligned} \tag{5.98}$$

that is,
$$\operatorname{tr}_\Sigma^- \boldsymbol{\eta} - \eta_\Sigma = 0, \tag{5.99}$$
$$\operatorname{tr}_\Sigma^+ \boldsymbol{\eta} - \alpha \eta_\Sigma = 0. \tag{5.100}$$

In particular, expressions (5.99) and (5.100) demonstrate that $\operatorname{tr}_\Sigma^\pm \boldsymbol{\eta} \in L^2(\Sigma)$, and that
$$\operatorname{tr}_\Sigma^+ \boldsymbol{\eta} - \alpha \operatorname{tr}_\Sigma^- \boldsymbol{\eta} = 0. \tag{5.101}$$

Properties (5.91), (5.94), and (5.101) imply that $a^*(\hat{\boldsymbol{\theta}}, \boldsymbol{\eta}) = 0 \ \forall \hat{\boldsymbol{\theta}} \in L^*$, and thus, by lemma 5.12, that $\boldsymbol{\eta} = 0$. Since $\boldsymbol{\eta}$ vanishes, also $\eta_\Sigma$ and $\boldsymbol{\eta}_s$ vanish, due to expressions (5.99) and (5.96), which finally proves the claim. $\square$

We then finally arrive at the well-posedness result for our variational formulation of initial–boundary-value problem (5.21).

**Theorem 5.16.** *With $a$ and $l$ as in definitions (5.17), where function $\alpha$ satisfies bound (5.4), with space $V$ as in definition (5.39), and with $L = L^2(Q)^{d+1} \times L^2(\Sigma) \times L^2(Q_0)^{d+1}$, the variational problem to find $\boldsymbol{\xi} \in V$ such that*
$$a(\hat{\boldsymbol{\eta}}, \boldsymbol{\xi}) = l(\hat{\boldsymbol{\eta}}) \qquad \forall \hat{\boldsymbol{\eta}} \in L \tag{5.102}$$
*has a unique solution satisfying*
$$\|\boldsymbol{\xi}\| \leq \frac{8\tau e}{1 - \alpha_M} \|l\|. \tag{5.103}$$

*Proof.* By the Cauchy–Schwarz inequality and trace lemma 5.11, $a$ and $l$ are continuous on $L \times V$ and $L$, respectively. Theorem 2.1 together with lemmas 5.14 and 5.15 then yields well-posedness of variational problem (5.102) and the bound (5.103). $\square$

**Acknowledgements.** The authors thank Rainer Picard for helpful discussions and for the inspiration that his work on evo-systems have had for the work presented here. We are also deeply thankful for the careful reading of the manuscript by the anonymous reviewers and their many constructive and helpful suggestions for improvements. Thanks also to Patrick Ciarlet for making us aware of Bourhrara's work [3] on the neutron transport equation. Funding for this work was partly provided by the Swedish Research Council, grant 2018-03546.

# Appendix

## The proof of lemma 5.6

**Lemma 5.6.** *Let $\boldsymbol{h} \in C^{0,\mu}(\partial Q)^n$ be a Hölder continuous function with exponent $\mu \in (1/2, 1]$. Then there is a constant $C$ such that for any $u \in H^{1/2}(\partial Q)$,*

$$\|\boldsymbol{h}u\|_{H^{1/2}(\partial Q)^n} \leq C \|u\|_{H^{1/2}(\partial Q)}. \tag{A.104}$$

*Proof.* Recalling the definition of the norm on $H^{1/2}(\partial Q)^n$ (note that $Q \subset \mathbb{R}^{d+1}$), we have that

$$\|\boldsymbol{h}u\|_{H^{1/2}(\partial Q)^n}^2 = \int_{\partial Q} |\boldsymbol{h}u|^2 + \int_{\partial Q \times \partial Q} \frac{|\boldsymbol{h}(x)u(x) - \boldsymbol{h}(y)u(y)|^2}{|x - y|^{d+1}}. \tag{A.105}$$

The inequality

$$\begin{aligned}|\boldsymbol{h}(x)u(x) - \boldsymbol{h}(y)u(y)|^2 &= |\boldsymbol{h}(x)u(x) - \boldsymbol{h}(y)u(x) + \boldsymbol{h}(y)u(x) - \boldsymbol{h}(y)u(y)|^2 \\ &\leq 2\left(|u(x)|^2 |\boldsymbol{h}(x) - \boldsymbol{h}(y)|^2 + |\boldsymbol{h}(y)|^2 |u(x) - u(y)|^2\right)\end{aligned} \tag{A.106}$$

implies that

$$\begin{aligned}\|\boldsymbol{h}u\|_{H^{1/2}(\partial Q)^n}^2 &\leq \int_{\partial Q} |\boldsymbol{h}u|^2 + 2 \int_{\partial Q \times \partial Q} \frac{|\boldsymbol{h}(y)|^2 |u(x) - u(y)|^2}{|x - y|^{d+1}} \\ &\quad + 2 \int_{\partial Q \times \partial Q} \frac{|u(x)|^2 |\boldsymbol{h}(x) - \boldsymbol{h}(y)|^2}{|x - y|^{d+1}}.\end{aligned} \tag{A.107}$$

Since $\partial Q$ is compact (closed and bounded) and $\boldsymbol{h}$ is (Hölder) continuous, there is a $\overline{h}$ such that

$$|\boldsymbol{h}| \leq \overline{h} \quad \text{on } \partial Q, \tag{A.108}$$

which substituted into inequality (A.107) yields that

$$\begin{aligned}\|\boldsymbol{h}u\|_{H^{1/2}(\partial Q)^n}^2 &\leq \overline{h}\|u\|_{H^{1/2}(\partial Q)}^2 + \overline{h} \int_{\partial Q \times \partial Q} \frac{|u(x) - u(y)|^2}{|x - y|^{d+1}} \\ &\quad + 2 \int_{\partial Q \times \partial Q} \frac{|u(x)|^2 |\boldsymbol{h}(x) - \boldsymbol{h}(y)|^2}{|x - y|^{d+1}}.\end{aligned} \tag{A.109}$$

Note that, if $I : \partial Q \to \overline{\mathbb{R}}$, defined by

$$I(x) = 2 \int_{\partial Q} \frac{|\boldsymbol{h}(x) - \boldsymbol{h}(y)|^2}{|x - y|^{d+1}} \, \mathrm{d}S_y, \tag{A.110}$$

is essentially bounded, the conclusion follows by Fubini's theorem.

By assumption, there is a constant $C_h$ and some $\mu \in (1/2, 1]$ such that

$$|\boldsymbol{h}(x) - \boldsymbol{h}(y)| \leq C_h |x - y|^\mu \quad \text{for all } x, y \in \partial Q. \tag{A.111}$$

Thus,

$$I(x) \leq 2C_h^2 \int_{\partial Q} |x - y|^{2\mu - d - 1} \, \mathrm{d}S_y. \tag{A.112}$$

By assumption, $\partial Q \subset \mathbb{R}^{d+1}$ is bounded, so

$$|x - y| \leq \operatorname{diam} \partial Q < \infty, \tag{A.113}$$

for any $x, y \in \partial Q$.

We will estimate the integral on the right side of inequality (A.112) by dividing $\partial Q$ into two parts, inside and outside a $\delta$ neighborhood of $x$, respectively. For this, since $Q$ lies locally on one side of its Lipschitz continuous boundary $\partial Q$, there are $\delta > 0$ and $L$ such that, for any $x \in \partial Q$, we have the bound

$$|x - y| \geq \delta \quad \text{for any } y \in \partial Q \setminus B_\delta(x), \tag{A.114}$$

where $B_\delta(x) \subset \mathbb{R}^{d+1}$ denotes the open ball with radius $\delta$ centered at $x$. Moreover, there is a local coordinate system and an $L$-Lipschitz continuous map $\theta : \omega \stackrel{\text{def}}{=} \{ \chi \in \mathbb{R}^d \mid |\chi| < 1/2 \} \to \mathbb{R}$ such that each point on $\partial Q \cap B_\delta(x)$ has a coordinate representation $(\chi, \theta(\chi))$ for some $\chi \in \omega$, and

$$\int_{\partial Q \cap B_\delta(x)} |x - y|^{2\mu-d-1} \, dS_y \leq \int_\omega |(\chi', \theta(\chi')) - (\chi, \theta(\chi))|^{2\mu-d-1} \sqrt{1 + \sum_{i=1}^d \left(\frac{\partial \theta}{\partial \chi_i}(\chi)\right)^2} \, d\chi \quad (A.115)$$

in which $(\chi', \theta(\chi'))$, $\chi' \in \omega$ is the coordinate representation of $x$ in the local coordinate system.

Now, considering first $y \in \partial Q \setminus B_\delta(x)$, by the bounds (A.113) and (A.114), it follows that

$$|x - y|^\lambda \leq \max\{\text{diam} \, \partial Q, 1/\delta\}^\lambda, \quad (A.116)$$

for any $\lambda \in \mathbb{R}$, which implies that

$$\int_{\partial Q \setminus B_\delta(x)} |x - y|^{2\mu-d-1} \, dS_y \leq |\partial Q| \max\{\text{diam} \, \partial Q, 1/\delta\}^{2\mu-d-1} < \infty, \quad (A.117)$$

where $|\partial Q|$ denotes the surface area (measure) of $\partial Q$.

In the case complementary to bound (A.117), we need to estimate the right side of inequality (A.115). Since $\theta$ is $L$-Lipschitz, it holds that

$$\sqrt{1 + \sum_{i=1}^d \left(\frac{\partial \theta}{\partial \chi_i}(\chi)\right)^2} \leq \sqrt{1 + dL^2}, \quad (A.118)$$

for almost any $\chi \in \omega$, and

$$|\chi' - \chi| \leq \sqrt{|\chi' - \chi|^2 + |\theta(\chi') - \theta(\chi)|^2} \leq \sqrt{1 + L^2}|\chi' - \chi|, \quad (A.119)$$

for any $\chi \in \omega$. Since $\sqrt{1 + L^2} \geq 1$, the bound (A.119) implies that

$$|(\chi', \theta(\chi')) - (\chi, \theta(\chi))|^\lambda = \sqrt{|\chi' - \chi|^2 + |\theta(\chi') - \theta(\chi)|^2}^\lambda \leq \sqrt{1 + L^2}^\lambda |\chi' - \chi|^\lambda, \quad (A.120)$$

for any $\lambda \in \mathbb{R}$ and $\chi \in \omega$. By introducing the bounds (A.118) and (A.120) into expression (A.115), we find that

$$\begin{aligned}
\int_{\partial Q \cap B_\delta(x)} |x - y|^{2\mu-d-1} \, dS_y &\leq C_{d,\mu,L} \int_\omega |\chi' - \chi|^{2\mu-d-1} \, d\chi \\
&\leq C_{d,\mu,L} \int_{|\chi'-\chi|<1} |\chi' - \chi|^{2\mu-d-1} \, d\chi \\
&= C_{d,\mu,L} \int_0^1 s^{2\mu-d-1} A_d s^{d-1} \, ds \quad (A.121) \\
&= C_{d,\mu,L} A_d \int_0^1 s^{2\mu-2} \, ds \\
&= C_{d,\mu,L} A_d \frac{1}{2\mu - 1} < \infty,
\end{aligned}$$

where $C_{d,\mu,L} = \sqrt{1 + dL^2}\sqrt{1 + L^2}^{2\mu-d-1}$, $A_d$ denotes the surface area of the $d$-dimensional Euclidean unit ball, we used that $\omega \subset \{\chi \in \mathbb{R}^d \mid |\chi'-\chi| < 1\}$ in the second inequality, and that $\mu \in (1/2, 1]$ in the last inequality.

Introducing the bounds (A.121) and (A.117) into expression (A.112), we find that

$$\begin{aligned}
I(x) &\leq 2C_h^2 \int_{\partial Q} |x - y|^{2\mu-d-1} \, dS_y \\
&= 2C_h^2 \int_{\partial Q \cap B_\delta(x)} |x - y|^{2\mu-d-1} \, dS_y + 2C_h^2 \int_{\partial Q \setminus B_\delta(x)} |x - y|^{2\mu-d-1} \, dS_y \quad (A.122) \\
&\leq 2C_h^2 \left(\frac{C_{d,\mu,L} A_d}{2\mu - 1} + |\partial Q| \max\{\text{diam} \, \partial Q, 1/\delta\}^{2\mu-d-1}\right) \stackrel{\text{def}}{=} C_{d,\mu,h,\partial Q} < \infty,
\end{aligned}$$

for any $\boldsymbol{x} \in \partial Q$; that is, $I$ is bounded on $\partial Q$. We may therefore apply Fubini's theorem to the second integral in inequality (A.109), invoke the bound (A.122) (recall definition (A.110)), and conclude that

$$
\begin{aligned}
\|\boldsymbol{h}u\|^2_{H^{1/2}(\partial Q)^n} &\leq \overline{h}\|u\|^2_{H^{1/2}(\partial Q)} + \overline{h}\int_{\partial Q \times \partial Q} \frac{|u(\boldsymbol{x}) - u(\boldsymbol{y})|^2}{|\boldsymbol{x} - \boldsymbol{y}|^{d+1}} \\
&\quad + 2\int_{\partial Q \times \partial Q} \frac{|u(\boldsymbol{x})|^2|\boldsymbol{h}(\boldsymbol{x}) - \boldsymbol{h}(\boldsymbol{y})|^2}{|\boldsymbol{x} - \boldsymbol{y}|^{d+1}} \\
&\leq \overline{h}\|u\|^2_{H^{1/2}(\partial Q)} + \overline{h}\int_{\partial Q \times \partial Q} \frac{|u(\boldsymbol{x}) - u(\boldsymbol{y})|^2}{|\boldsymbol{x} - \boldsymbol{y}|^{d+1}} \\
&\quad + C_{d,\mu,h,\partial Q}\int_{\partial Q} |u(\boldsymbol{x})|^2 \\
&\leq \left(\overline{h} + \max\{\overline{h}, C_{d,\mu,h,\partial Q}\}\right) \|u\|^2_{H^{1/2}(\partial Q)}.
\end{aligned}
\tag{A.123}
$$

□